%% file: gear_design.tex
\documentclass[12pt]{amsart}

\input{header_basic.tex}

\input{header_article.tex}
\input{header_subtle.tex}

\usepackage{enumitem}
\usepackage{color}
\usepackage{subfig, caption}
\usepackage{wrapfig}
\usepackage{pdflscape}
\usepackage{array}
\usepackage{MnSymbol}

\title[Mathematical overview and applications of gear design]{A mathematical overview and some applications of gear design}
\author{Elisabetta A. Matsumoto and Henry Segerman} 
\date{\today}

\begin{document}

\begin{abstract}
In this paper we give a brief overview of the geometry of involute gears, from a mathematical more than an engineering perspective. 
We also list some of the many variant geared mechanisms and discuss some of our 3D printed mechanisms.
\end{abstract}

\maketitle
\section{Introduction}

Gears are some of the most ubiquitous mechanical parts.
They transfer torque between axles, often applying leverage as they do so, trading angular speed for torque and vice versa.
Engineering texts on gears include many details that are used to characterize standardized commercial gears. In this paper, we aim to describe the construction of gears from a purely geometric perspective.

\begin{figure}[h]
\centering
\subfloat[]{
\includegraphics[width=0.3\textwidth]{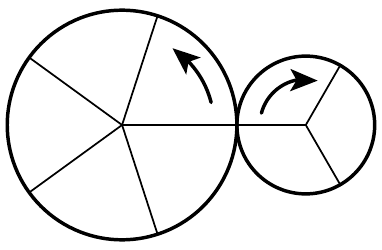}
}
\subfloat[]{
\includegraphics[width=0.3\textwidth]{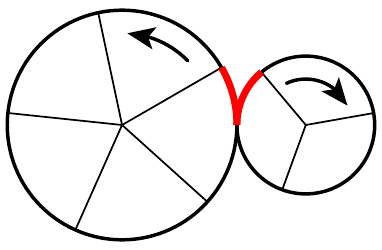}
}
\subfloat[]{
\includegraphics[width=0.3\textwidth]{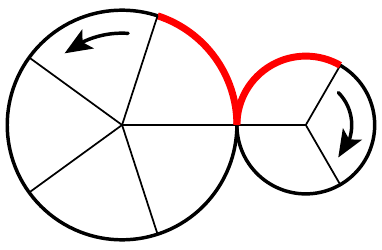}
}
\caption{Rolling circles. In each diagram the red curves have equal length.}
\label{Fig:DrivenDriving}
\end{figure}

Consider two circles of radius $r_1$ and $r_2$ that rotate without slipping against each other, as shown in \reffig{DrivenDriving}. The fundamental principle of gears is that as the gears (or, for now, circles) rotate, the arc length swept out by the driving gear (with radius $r_1$, say) must be equal to the arc length swept out by the driven gear (with radius $r_2$). For a circle of radius $r$ that has rotated by angle $\theta$ radians, the arc length travelled is $s=r\theta$. 
Thus the angle that the driven gear rotates is 
\[\theta_2=(r_1/r_2) \theta_1\] 

\begin{figure}[h]
\centering
\labellist
\small\hair 2pt
\pinlabel $R$ [l] at 56 32
\endlabellist
\includegraphics[height=6cm]{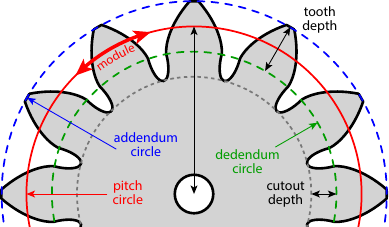}
\caption{Anatomy of a gear.}
\label{Fig:anatomy}
\end{figure}

We list some engineering terms that describe different parts of a gear, shown in \reffig{anatomy}. These terms are not universally adopted, and other sources may vary. A pair of gears is based on two curves rotating around fixed axles that roll against each other without slipping.  Such a curve, shown in red in the figure, is called the \emph{pitch curve}. For now, we will consider circular gears, but in Sections \ref{Sec:AcircularPitchCurves} through \ref{Sec:AlienGears} we will look at gears with acircular pitch curves. The rest of the terminology refers to the specifics of the gear teeth. The teeth are inscribed in a curve called the \emph{addendum curve}, shown in blue. The bases of the teeth lie along the \emph{dedendum curve} (sometimes called a \emph{base circle}), shown in green. The distance between addendum and dedendum curves is the \emph{tooth depth}. Some gears have an additional cutout between teeth for clearance, whose depth is called the \emph{cutout depth}. The most important measurement is the \emph{module} $m = R\frac{n}{2\pi}$, the fraction of the pitch curve of radius $R$ from a gear with $n$ teeth for each tooth.

\subsection*{Acknowledgements} This material is based in part upon work supported by the National Science Foundation under Grant No. DMS-1439786 and the Alfred P. Sloan Foundation award G-2019-11406 while the authors were in residence at the Institute for Computational and Experimental Research in Mathematics in Providence, RI, during the Illustrating Mathematics program.
In addition, the first author was supported in part by National Science Foundation grant DMR-1847172 and a Cottrell Scholars Award from the Research Corporation for Science Advancement. 
The second author was supported in part by National Science Foundation grant DMS-2203993.

\section{Involute gears}
\label{Sec:Involute}

\subsection{Adding teeth to a gear}

\begin{wrapfigure}[33]{r}{0.43\textwidth}
\vspace{-22pt}
\centering
\subfloat[Line of force (green) transmitted by gear teeth at pressure angle $\alpha$.]{
\label{Fig:PressureAngle}
\labellist
\small\hair 2pt
\pinlabel $\alpha$ [l] at 110 35
\endlabellist
\includegraphics[width=0.3\textwidth]{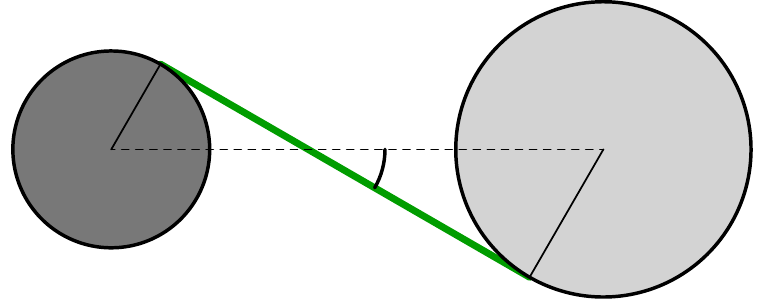}
}

\subfloat[Transmitting force along a line.]{
\label{Fig:InvolutePair}
\includegraphics[width=0.3\textwidth]{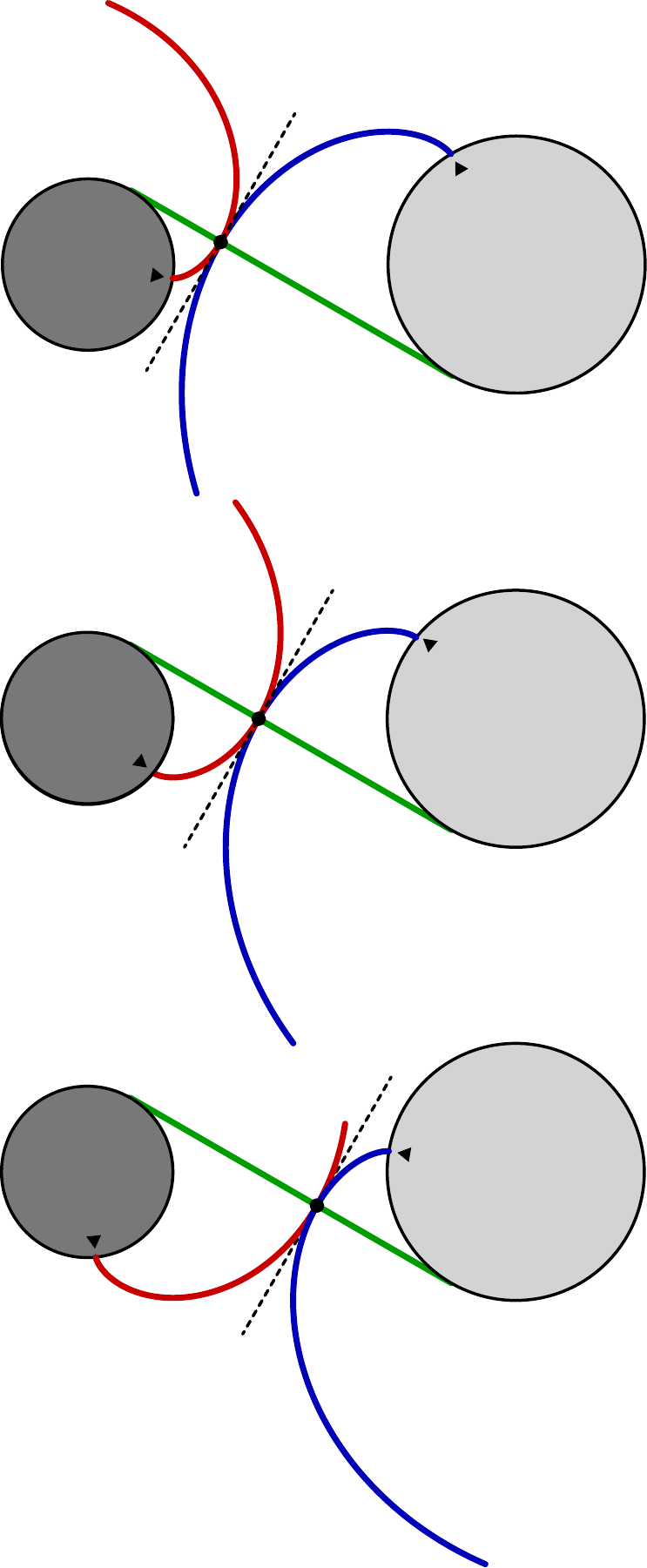}
}
\caption{A pair of involutes based on circles of different radii rotating against each other. The dot shows the point of contact between the two curves. }
\label{Fig:force}
\end{wrapfigure}

There are two main concerns when adding teeth to a gear. First we need to make sure that the teeth mesh together, and second the gears should efficiently transfer torque from the driving gear to the driven gear. Teeth meshing makes the two gears rotate in registry -- so that their respective pitch curves act as if they rotate without slipping. To ensure that the teeth in our gears mesh, we must have the same number of teeth per unit arc length on the two gears. So if our driving gear has $n_1$ teeth, our driven gear must have $n_2=(r_1/r_2) n_1$ teeth. Of course both $n_1$ and $n_2$ must be integers. Frequently, gears that have not been precisely machined will use trapezoidal shaped teeth. This is practical for ensuring that teeth mesh, but the teeth won't contact each other smoothly, resulting in inefficient transmission of force.

\subsection{Involutes}
\begin{figure}[h]
\centering
\includegraphics[height=4cm]{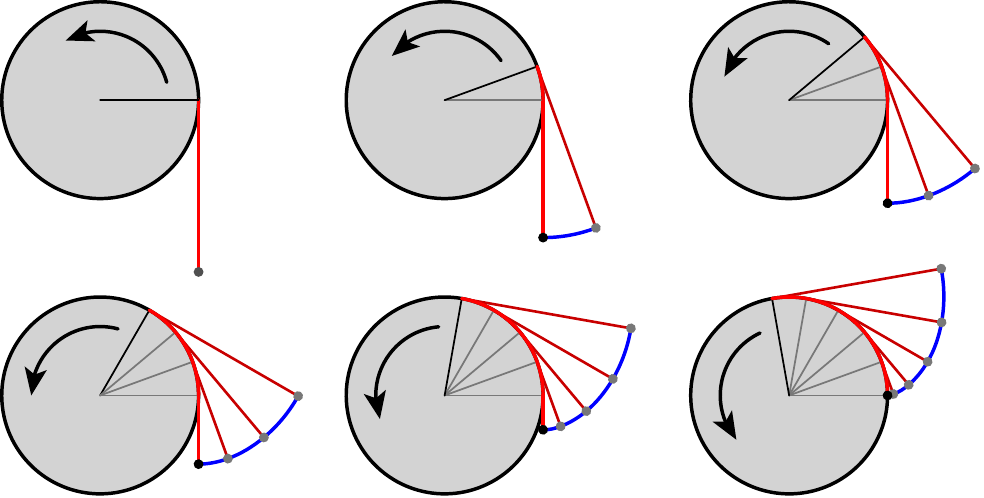}
\caption{Construction of an involute curve.}
\label{Fig:involute_construction}
\end{figure}

Almost all circular gears use \emph{involutes} for the tooth flanks. (See \reffig{involute_construction} for the construction. For an excellent animated introduction to involute gearing, see~\cite{Bencsik}.)
In physical gears this is critical for efficiency of torque transmission between components: throughout the rotation, we want equal torque to be transmitted from the driving gear to the driven gear. 

For any pair of circles of radii $r_1$ and $r_2$ which form the base (or dedendum) circles for a gear pair, there exists a line tangent to both along which force is transmitted. In \reffig{InvolutePair} this line is drawn in green. 
As the red curve rotates clockwise, it applies force to the blue curve along the green line, causing it to rotate anticlockwise.
As the gears rotate, the point of contact stays on the green line.
Moreover, the shared line of tangency between the curves (dashed line) is always perpendicular to the line of contact.
Together, these guarantee that the torque 
remains constant throughout the entire rotation of the gear. 

The angle $\alpha$ that the green line makes with the line connecting the axles of the gears is called the \emph{pressure angle}. See \reffig{PressureAngle}.
(Most commercial gears use a pressure angle of $20^\circ$, for its balance between force transmission and tooth strength.) 

The red and blue curves of \reffig{InvolutePair} are involutes.
\reffig{involute_construction} shows a construction of an involute: we can think of it as the curve traced out by the end of a piece of string as it is wrapped around another curve $\gamma(t)$. 
The idea is to take a line of length $\ell$ that is tangent to $\gamma$ at $t=t_0$ and draw a point. Then move (very slightly) down the curve $\gamma$ by arc length $a$ (at point $t=t_1$). Draw a point at the end of the line of length $\ell-a$ along the tangent at the point $\gamma(t_1)$. Repeat this process and connect the points. The result is an involute. 
If the curve $\gamma$ is a circle parametrized by $\gamma(\theta) = \{r \cos \theta, r \sin\theta\}$ then the shape of its involute is given by 
\begin{equation}
\boldsymbol{c}(\theta) = \{r(\cos \theta + \theta \sin \theta), r(\sin \theta - \theta \cos\theta)\}.
\end{equation}
More generally, the equation for an involute $c(t)$ to a curve $\gamma(t)$ at the parameter $t=t_0$ is:
\begin{equation}
\boldsymbol{c}_{t_0}(t) = \gamma(t) -\frac{\boldsymbol{c}(t)}{\vert \boldsymbol{c}(t) \vert} \int_{t_0}^{t} \left| \boldsymbol{c}(s)\right| ds.
\end{equation}

\subsection{Gear trains: module and pressure angle}
\label{Sec:teeth}

\begin{figure}[h]
\centering
\labellist
\small\hair 2pt
\pinlabel $\alpha$ at 70.2 46
\endlabellist
\includegraphics[height=4cm]{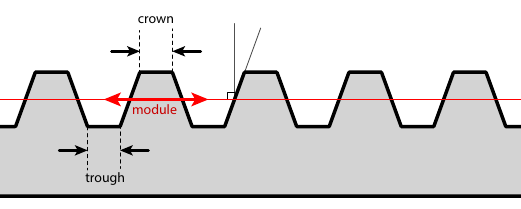}
\caption{Schematic of a rack with pressure angle $\alpha$.}
\label{Fig:RackSchematic}
\end{figure}

The genius behind involute teeth lies not only in their ability to transmit force, but in the manufacturing process that machines them. In effect, gears are cut by other gears. The die used to cut the teeth on a gear is called a \emph{rack}. This is a ``straightened out gear''.
A rack can be seen as the limit of an involute gear as the radius of the pitch circle approaches infinity. 
The geometry of a rack is defined by the module $m$ and the pressure angle $\alpha$. See \reffig{RackSchematic}. 

\begin{figure}[h]
\centering
\subfloat[Rack pair]{
\includegraphics[height=2.5cm]{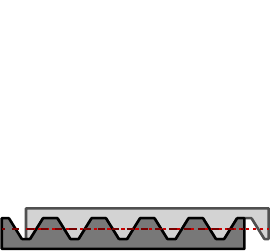}
}
\subfloat[Carving gears]{
\includegraphics[height=2.5cm]{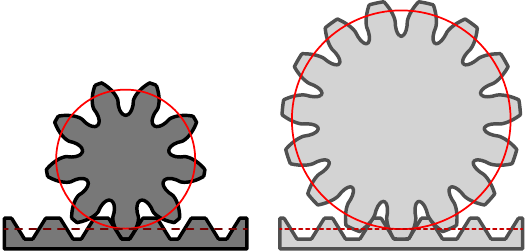}
}
\subfloat[Gear pair]{
\includegraphics[height=2.5cm]{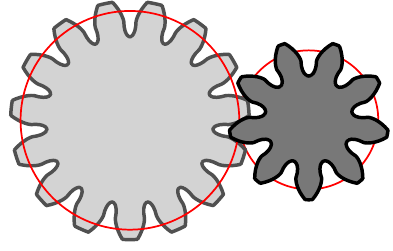}
}
\caption{Creating a pair of matching gears.}
\label{Fig:gear_pair}
\end{figure}

In general, a pair of racks that are complementary (\reffig{gear_pair}a), can be used to cut a pair of gears (\reffig{gear_pair}b) that will mesh together (\reffig{gear_pair}c). In practice, gears are cut by racks that are self-complementary, where the tooth half-width is half the module, meaning that the crown and trough lengths of the rack are identical. Thus as long as two gears have pitch circles whose circumferences are an integer number of times the module of the generating die, then those two gears will mesh with each other. This technique can be used to make gear trains or other networks of multiple gears of different sizes.

\section{Traditional gear variants}
\label{Sec:Variants}

Varying the construction given in \refsec{Involute} produces some closely related mechanisms.
These variant contexts often require subtly different gear flank shapes to maintain optimal torque transfer. 
Giving parametric descriptions of these shapes can be difficult.
In practice, engineers cut gears from dies to precisely match the complementary geometry, and therefore don't need parametric descriptions.
We will not worry too much about these details and here just give a catalogue of different commonly used gear types.

\subsection{Spur, helical, and herringbone gears}

\begin{figure}[h]
\centering
\includegraphics[width=\textwidth]{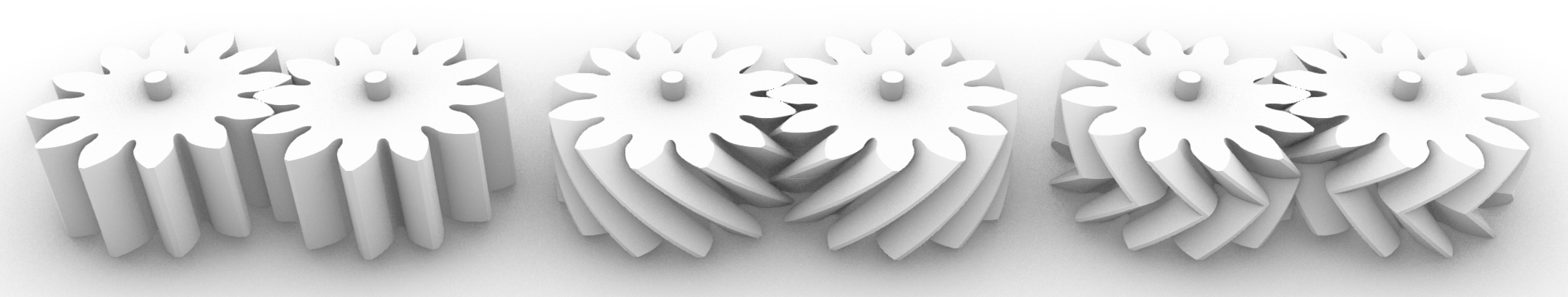}
\caption{From left to right: spur, helical, and herringbone gears.}
\label{Fig:SpurHelicalHerringbone}
\end{figure}

So far we have discussed essentially two-dimensional mechanisms, living in a plane. 
In real-world applications, the gears must be thickened up. The simplest way to thicken a planar gear is to extrude it perpendicular to the plane of the gear. (Mathematically, one would say that the three-dimensional shape is a product with an interval perpendicular to the plane.) This produces \emph{spur gears}. 
If we twist the planar gear as we extrude then we get \emph{helical gears}. 
If we glue two helical gears of opposite handedness together we get a \emph{herringbone gear}. Herringbone gears have the advantage that a meshing pair tends to stay in position. See \reffig{SpurHelicalHerringbone}.

\subsection{Bevel gears}

\begin{figure}[!h]
\centering
\subfloat[]{
\includegraphics[height = 5cm]{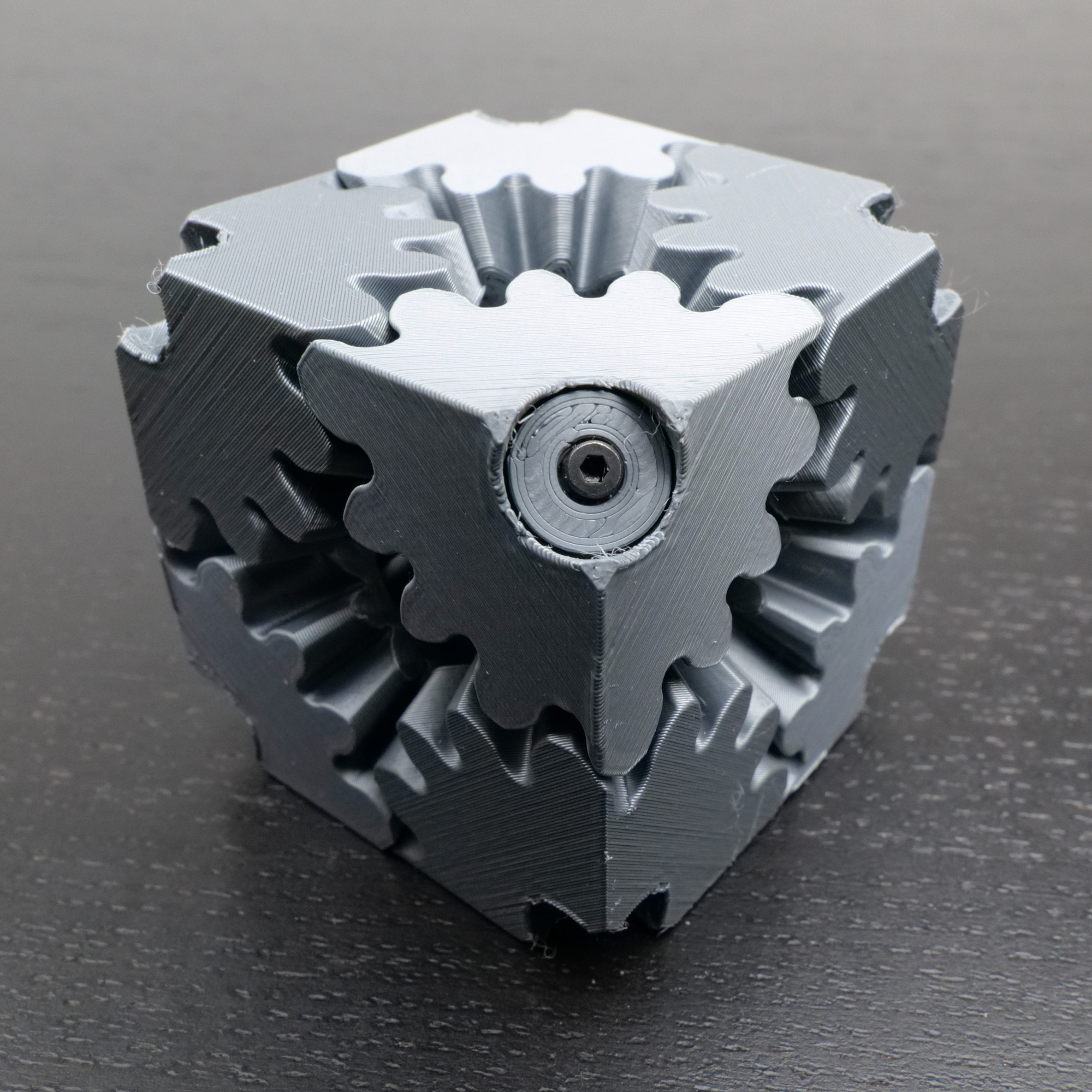}
}
\subfloat[]{
\includegraphics[height = 5cm]{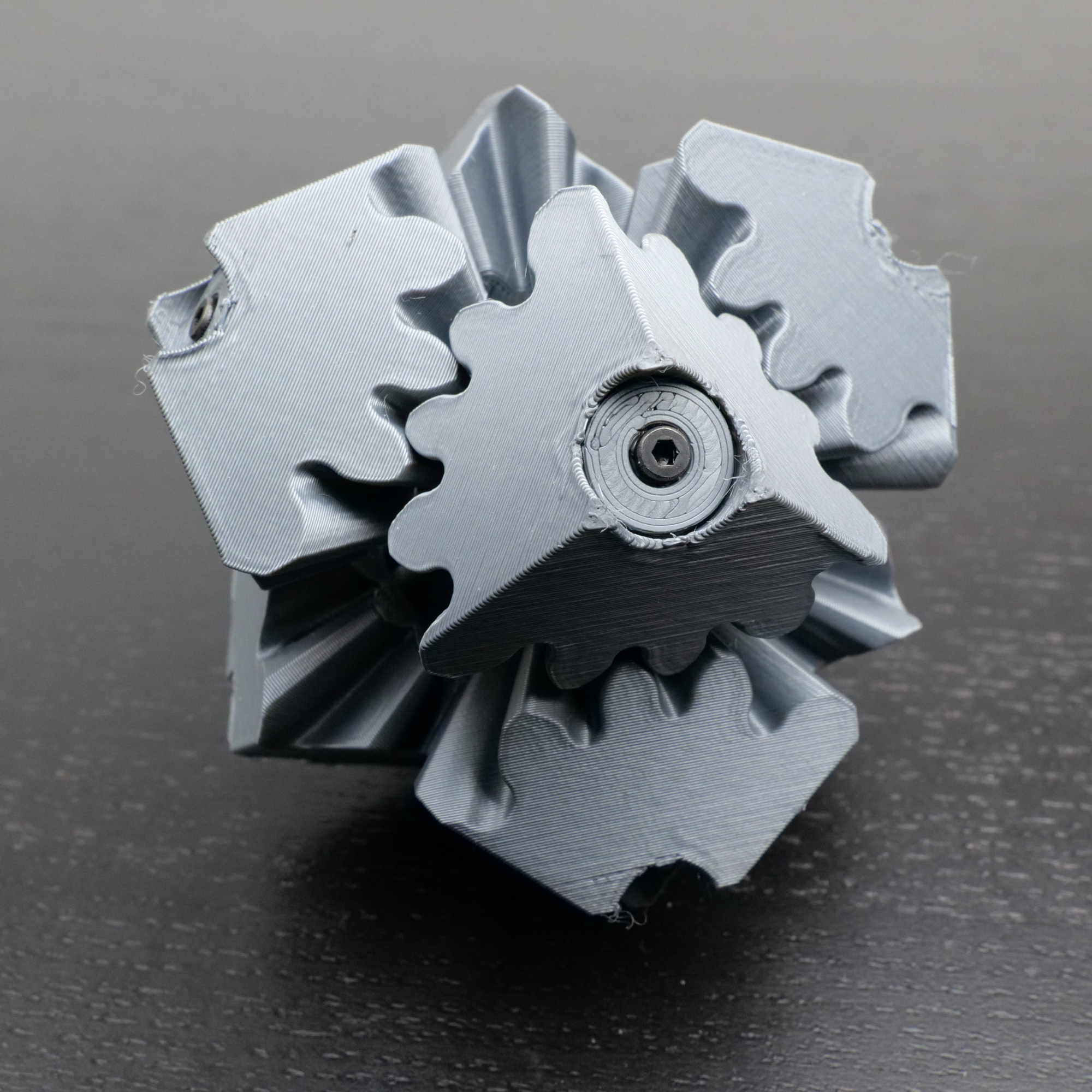}
}
\caption{Gear cube by the authors.}
\label{Fig:BrainGearCube}
\end{figure}

\begin{figure}[!h]
\centering
\subfloat[]{
\includegraphics[height = 5cm]{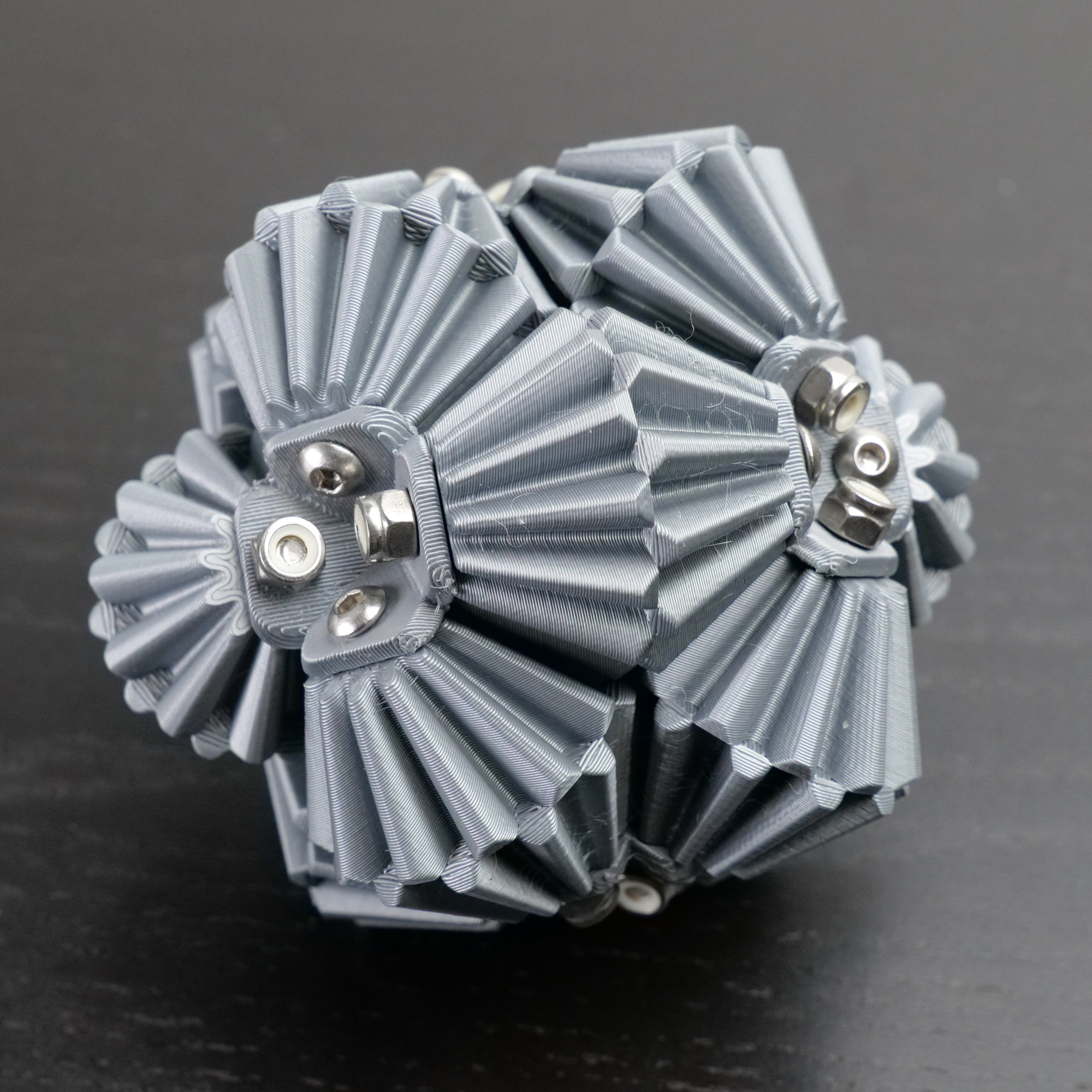}
}
\subfloat[]{
\includegraphics[height = 5cm]{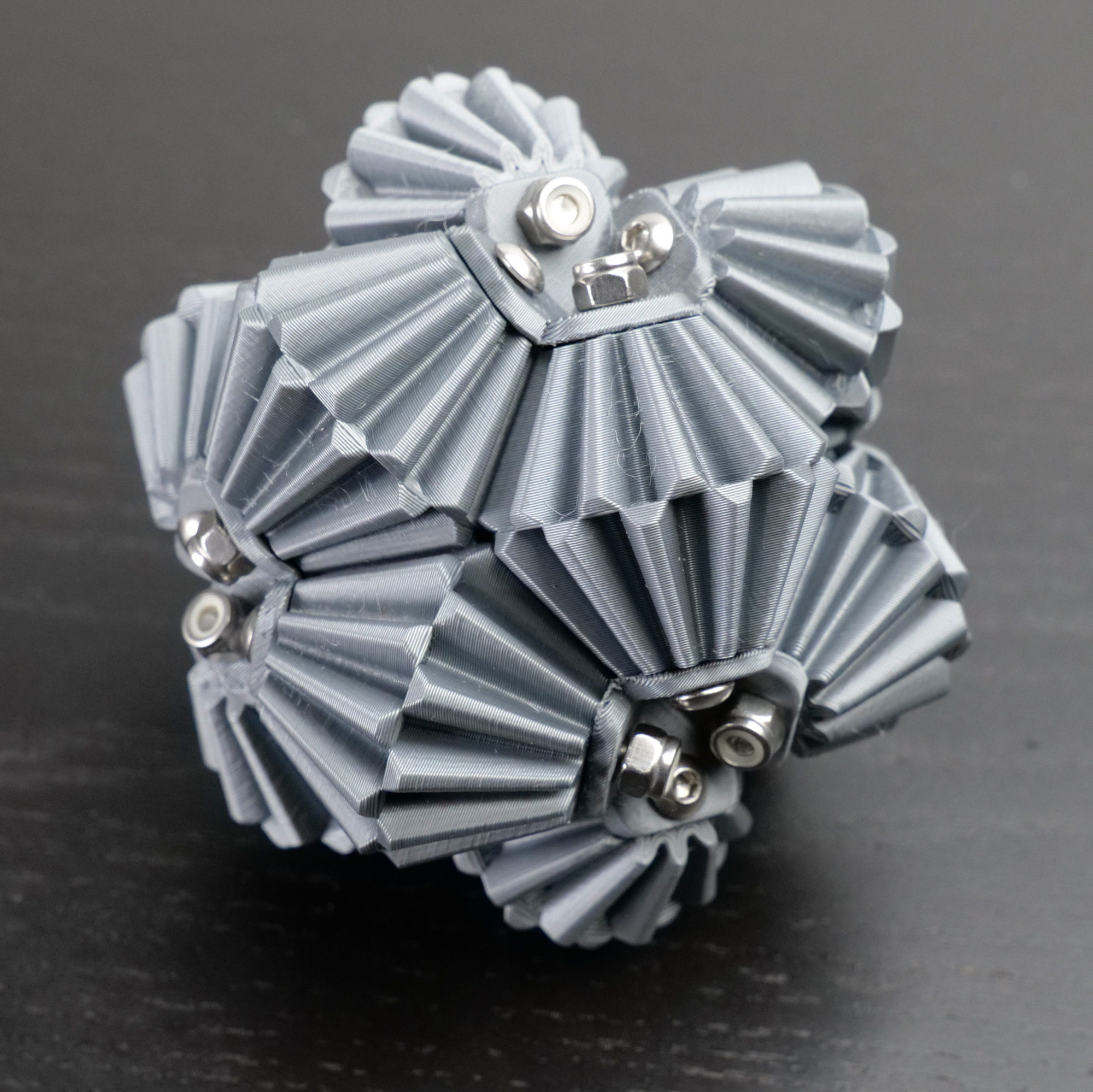}
}
\caption{Brain gear by the authors.}
\label{Fig:BrainGearEdges}
\end{figure}

Spur, helical, and herringbone gears rotate around parallel axles that are perpendicular to the plane in which the gears are constructed.
Moving the axes so that they meet at a point requires that the gears become \emph{bevel gears}.
Although not mathematically perfect, a good approximation to the correct shape for such a bevel gear can be formed by coning a planar involute gear to a point. (The correct shape seems to be a \emph{spherical involute}, which is constructed analogously to a planar involute~\cite{ParkLee}.)
In our experience this approximation is good enough for 3D printed mechanisms.
Figures \ref{Fig:BrainGearCube} and \ref{Fig:BrainGearEdges} show two designs based on bevel gears (see~\cite{BevelGearsVideo} for a video). 
In \reffig{BrainGearCube} the axles all meet at a single point at the center of a cube. 
In \reffig{BrainGearEdges}, the axles meet at the the vertices of an octahedron. 

\begin{figure}[!h]
\centering
\includegraphics[height = 9cm]{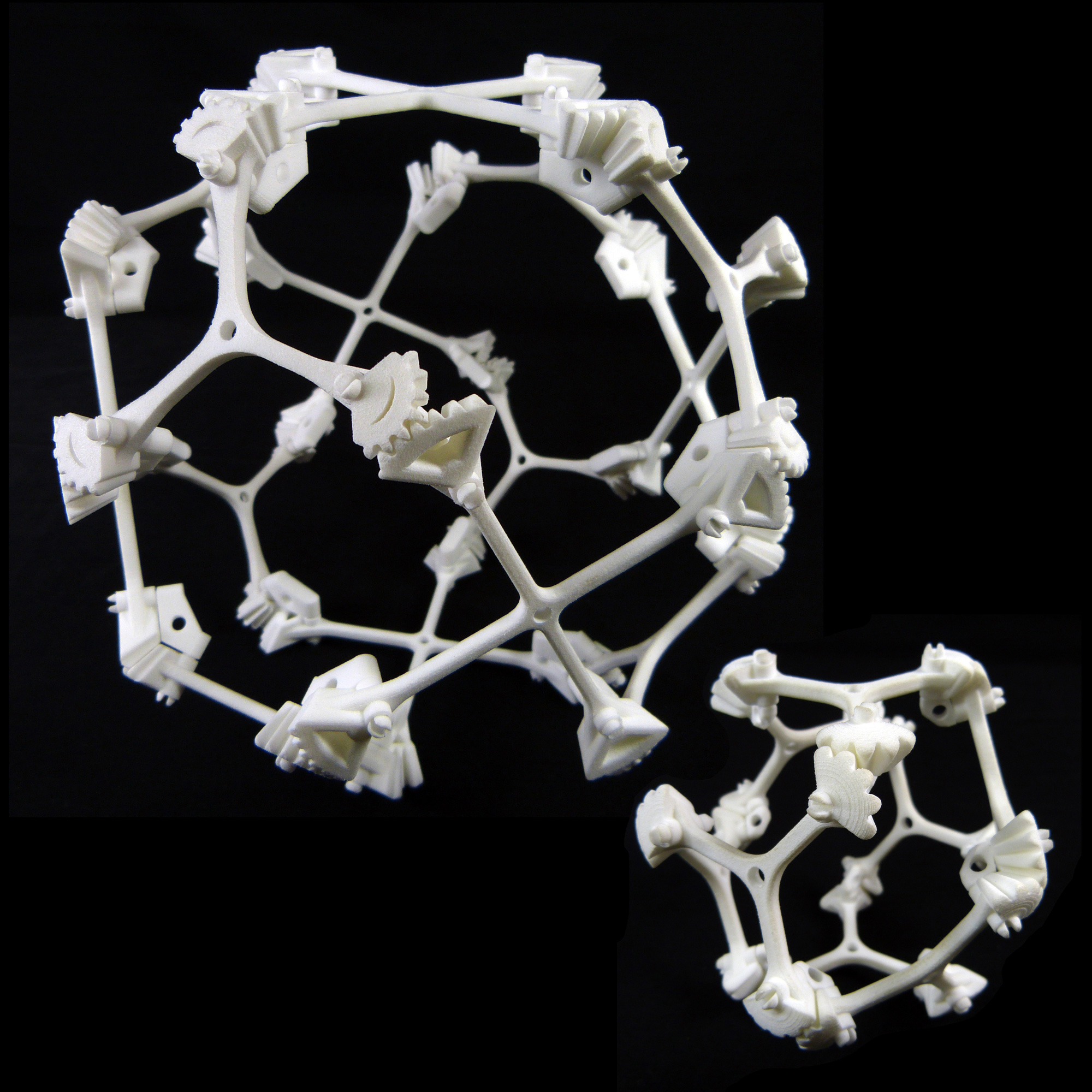}
\caption{Geared jitterbugs by the authors. Videos on these designs:~\cite{GearedJitterbugVideo, GearedCuboctahedralJitterbugVideo}}
\label{Fig:GearedJitterbugs}
\end{figure}

We have also used bevel gears to stabilize Buckminster Fuller's jitterbug mechanism, as shown in \reffig{GearedJitterbugs}.
For more details, see~\cite{GearedJitterbugs}.

\subsection{Skew and Worm gears}

\begin{figure}[!h]
\centering
\includegraphics[width=0.6\textwidth]{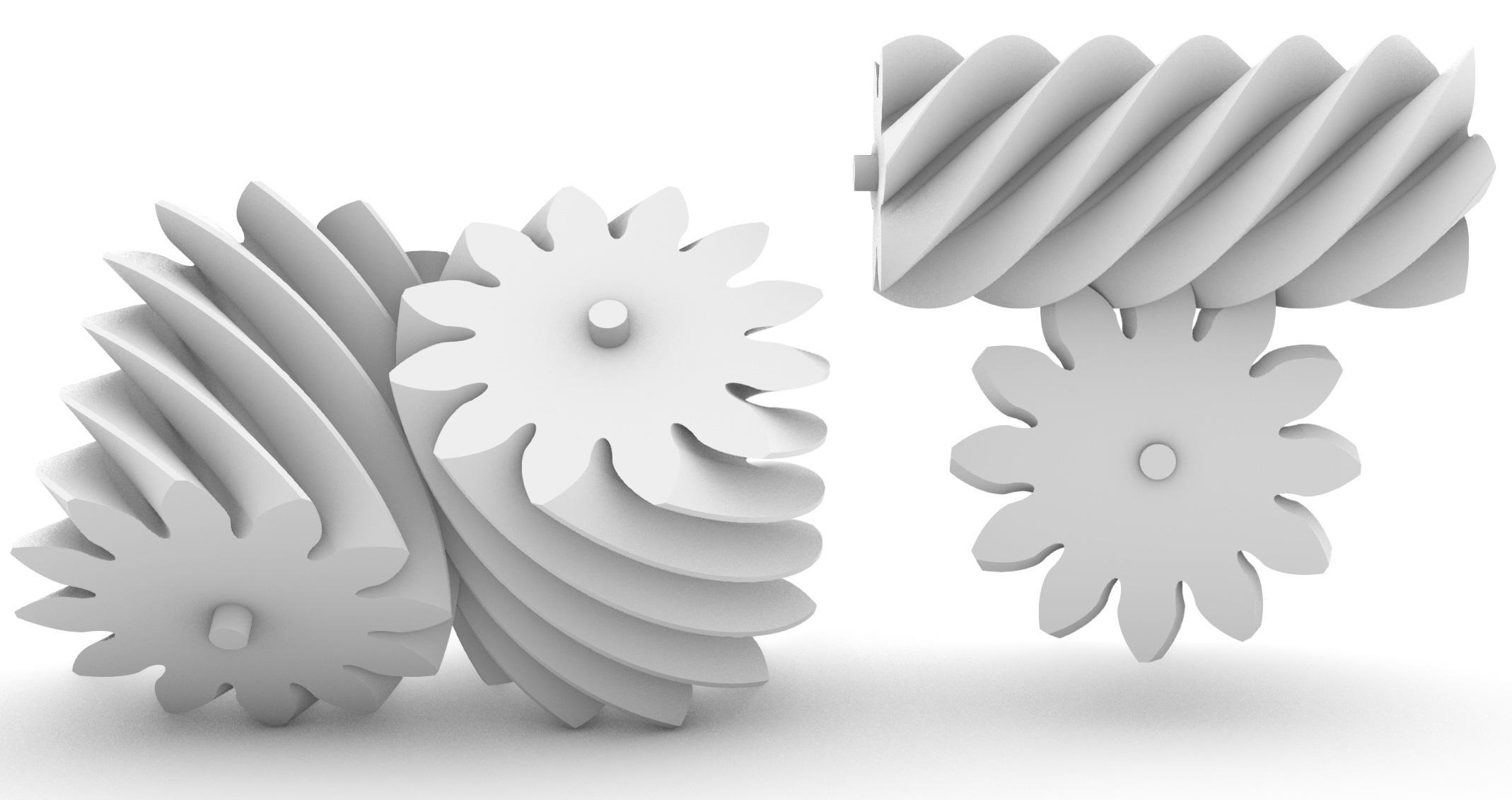}
\caption{Skew gears (left) and a worm drive (right).}
\label{Fig:SkewWorm}
\end{figure}

When the axles of a pair of meshing gears are neither parallel nor meet in a point, the gears become \emph{skew gears}. 
A common mechanism of this kind is a \emph{worm drive}, consisting of a worm (a skew gear similar to a screw) driving a \emph{worm wheel}. This gives a large reduction ratio. 
See \reffig{SkewWorm}.

\section{Acircular gears}

Acircular (or non-circular) gears are used in many engineering applications. Instead of having the pitch curve be a constant distance from the axle, acircular gears have a pitch curve $r(\theta)$ whose radius varies with the angle around the axle. The purpose of these gears is not to transmit torque with optimal efficiency but instead to control the gear ratio (as in a continuous transmission) or to control an additional motion (for example oscillations in a cam)~\cite{acirc_involute}.

\subsection{Acircular pitch curves}
\label{Sec:AcircularPitchCurves}

\begin{figure}[!h]
\centering
\labellist
\small\hair 2pt
\pinlabel $r_1(\theta_1)$ [l] at 606 72
\small\hair 2pt
\pinlabel $r_2=a-r_1(\theta_1)$ [l] at 56 72
\small\hair 2pt
\pinlabel $\theta_1$ [l] at 676 116
\small\hair 2pt
\pinlabel $\theta_2$ [l] at 106 116
\small\hair 2pt
\pinlabel $s_2(\theta_2)=s_1(\theta_1)$ [l] at 303 150
\endlabellist
\includegraphics[width = .9\textwidth]{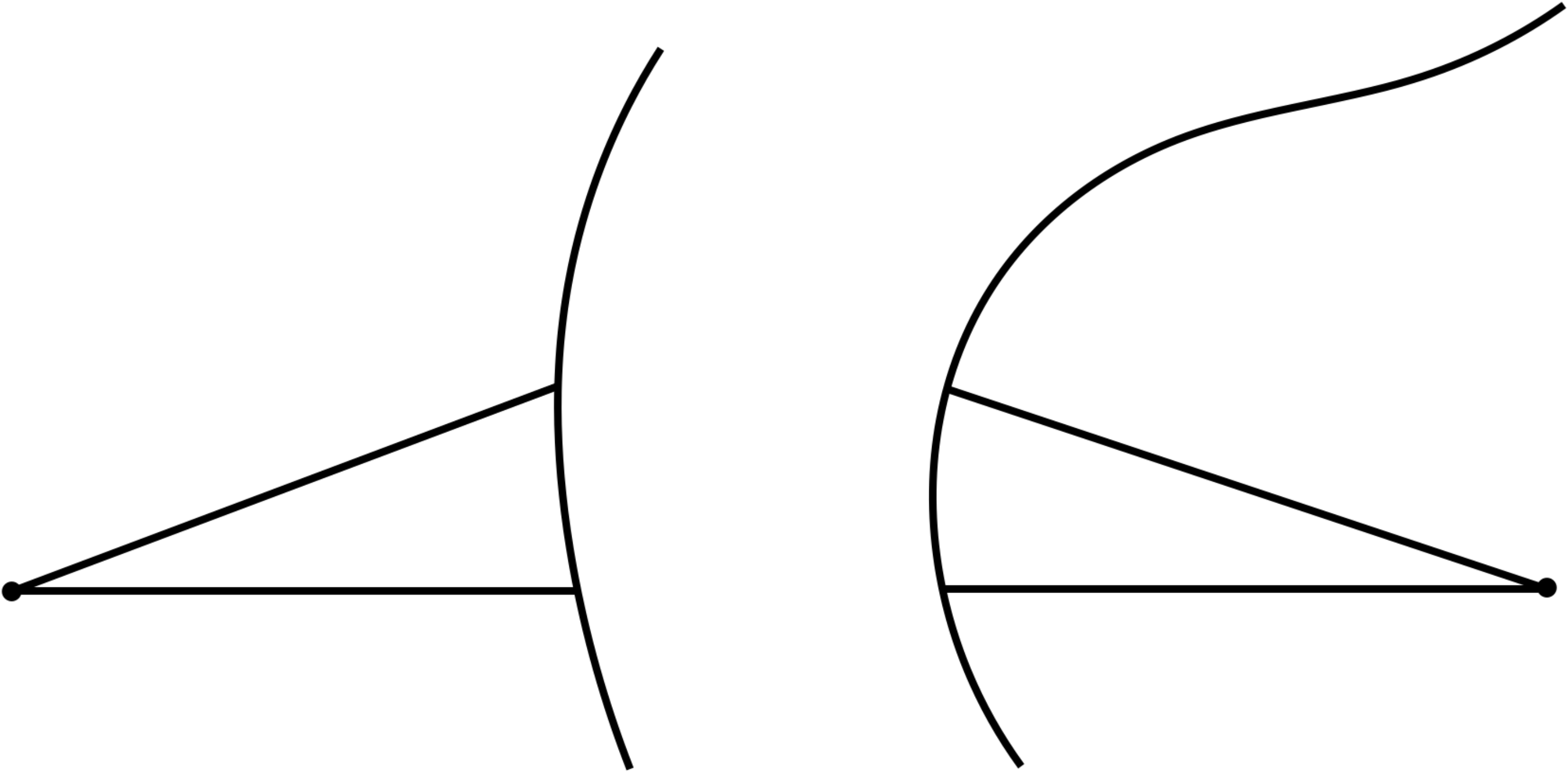}
\caption{Geometry of an acircular gear pair.}
\label{Fig:acirculara}
\end{figure}

We consider axles which are separated by a fixed distance, $a = r_1 + r_2$. 
If we know the pitch curve of the driving gear $r_1(\theta_1)$ then we can express $r_2$ in terms of $\theta_1$. However, to make the driven gear, we need to express $r_2$ as a function of its own rotation coordinate, $\theta_2$.
Since the arclength swept out by each gear in time $t$ must be the same, we end up with a relationship between their radii and rotation rates, as illustrated in \reffig{acirculara}: 
\begin{equation}
\label{Eqn:acirc_rate}
r_1 d\theta_1 = r_2 d\theta_2 = \big(a-r_1(\theta_1)\big)d\theta_2
\end{equation}
If we assume that the driving gear with radius $r_1(\theta_1)$ rotates at a constant rate, then we can integrate \refeqn{acirc_rate} to obtain the rotation rate of the driven gear:
 \[
 \theta_2+c = \int \frac{r_1(\theta_1)}{a-r_1(\theta_1)}d\theta_1.
 \] 
 This then enables us to calculate the pitch curve of the driven gear $r_2(\theta_2)$.

One classical example of an acircular gear pair is a nautilus gear pair, where the pitch curves are given by involutes of circles, as shown in \reffig{InvolutePair}.

\subsection{Variable speed gears}

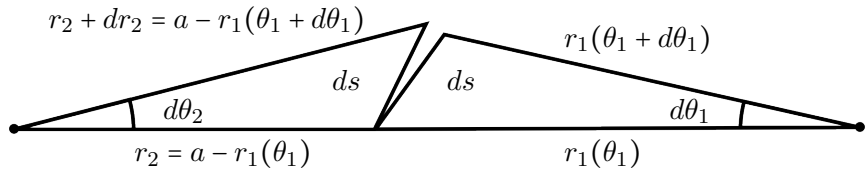
\begin{figure}[!h]
\labellist
\small\hair 2pt
\pinlabel $r_1(\theta_1)$ [l] at 576 32
\small\hair 2pt
\pinlabel $r_2+dr_2=a-r_1(\theta_1+d\theta_1)$ [l] at 36 172
\small\hair 2pt
\pinlabel $r_1(\theta_1+d\theta_1)$ [l] at 576 152
\small\hair 2pt
\pinlabel $r_2=a-r_1(\theta_1)$ [l] at 126 32
\small\hair 2pt
\pinlabel $d\theta_1$ [l] at 686 77
\small\hair 2pt
\pinlabel $d\theta_2$ [l] at 156 77
\small\hair 2pt
\pinlabel $ds$ [l] at 453 110
\small\hair 2pt
\pinlabel $ds$ [l] at 333 110
\endlabellist
\includegraphics[width = .9\textwidth]{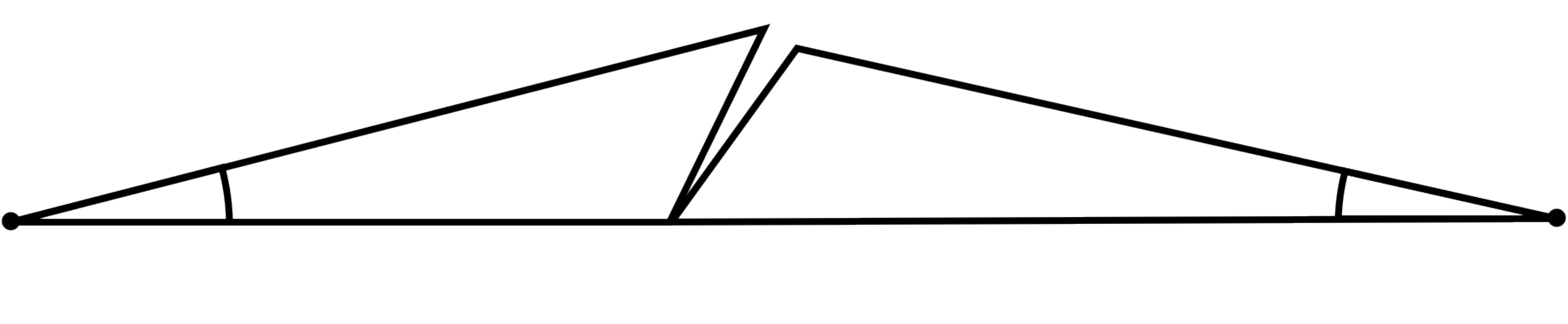}
\caption{Differential elements of an acircular gear pair.}
\label{Fig:acircularb}
\end{figure}

For certain applications, we might not know the shape of the pitch curve of our acircular gears, but we do know the desired (variable) relative rotation rates between the two gears.
That is, we require that the driven gear must rotate by some prescribed angle $\theta_2(\theta_1)$.
An infinitesimal interpretation of this condition is given in \reffig{acircularb}. Integrating this, we obtain
\[
\int_0^{\theta_1}\sqrt{r_1(s)^2+r_1'(s)^2}ds = \int_0^{\theta_2}\sqrt{r_2(t)^2+r_2'(t)^2}dt.
\]
This is unlikely to be analytically solvable. We give an iterative solution in \cite{GearedJitterbugs}.

We needed to use variable speed gearing in the design of the larger jitterbug in \reffig{GearedJitterbugs}. This mechanism is based on a cuboctahedral shape. Each pair of bevel gears consists of a gear at a corner of a square face of the cuboctahedron meeting a gear at a corner of a triangular face. To get to the open state the square rotates by 90 degrees while the triangle rotates by 120 degrees. However, the rotation speeds of the two gears are not constant during this motion. For more details, again see \cite{GearedJitterbugs}.

Note that both this construction and the construction given in \refsec{AcircularPitchCurves} take advantage of the fact that rolling without slipping means that the arclength swept out by the two gears in a given time must be the same.

\subsection{Cutting teeth on acircular gears}
\label{Sec:AcircularCutting}

\begin{figure}[h]
\centering
\subfloat[The rack must be tangent to both pitch curves during cutting.]{
\label{Fig:TangentRack}
\includegraphics[height=3.4 cm]{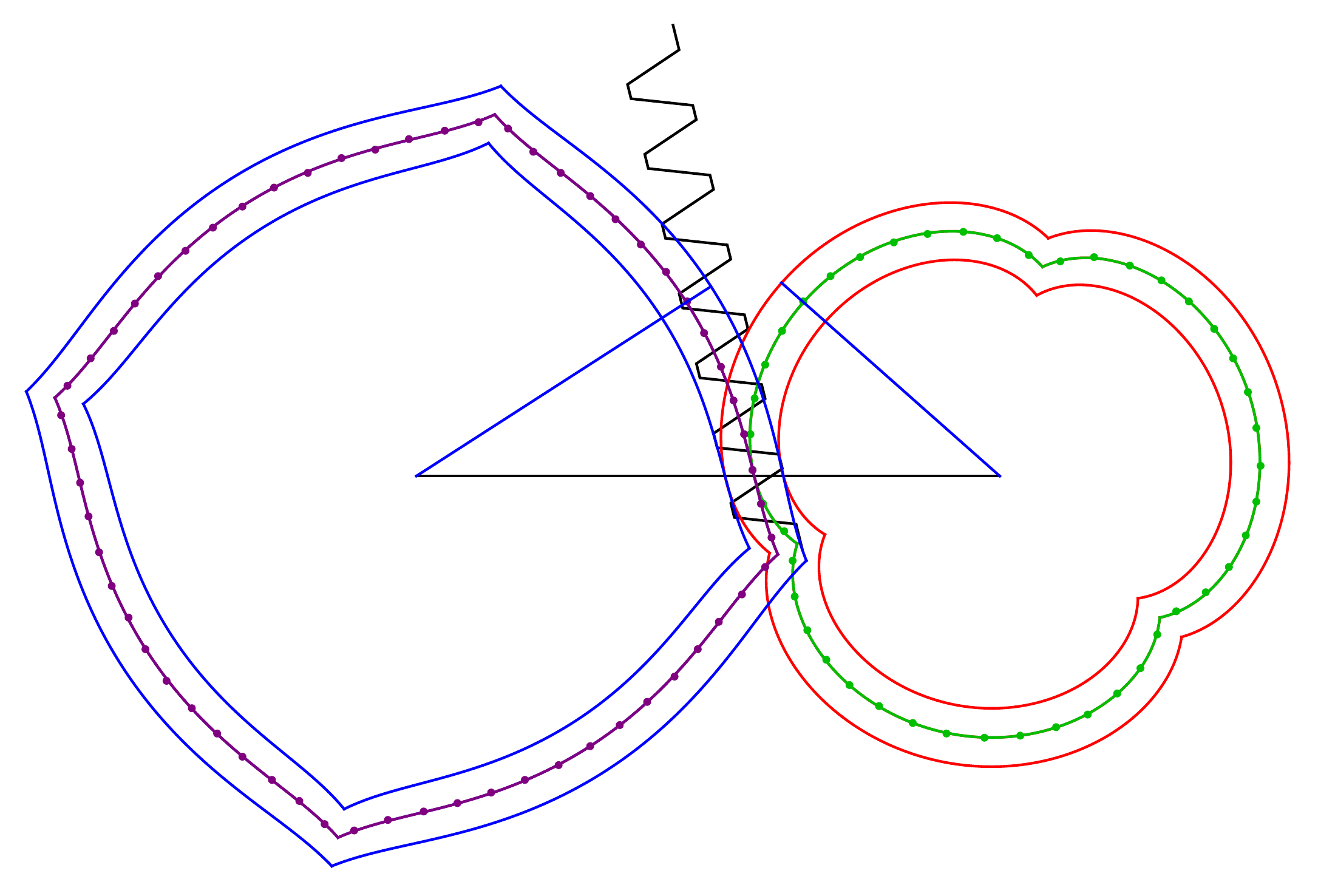}
}\hspace{5pt}
\subfloat[The tooth shapes are given by the envelope of the rack as it sweeps along each gear.]{
\label{Fig:Envelope}
\includegraphics[height=3.4 cm]{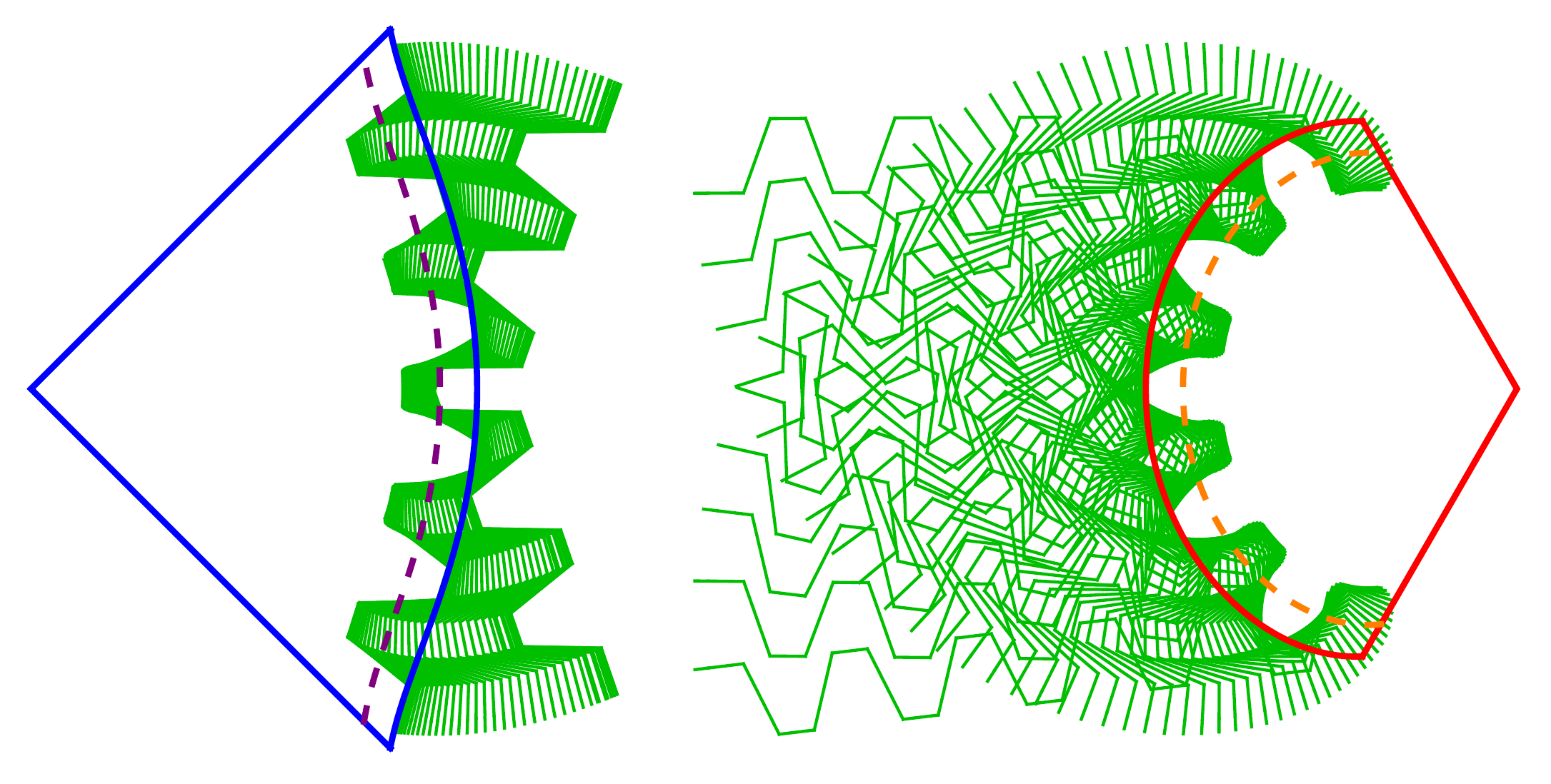}
}
\caption{Cutting teeth of acircular gears with a rack.}
\label{Fig:CuttingAcircular}
\end{figure}

The process of adding teeth to acircular gears is quite similar to that discussed in \refsec{teeth}. 
Again, the profile of the teeth is cut by a rack. 
The distance the rack travels as the gear rotates must be the same as the arclength swept out by the pitch curve during the motion. Additionally, the pitch curve in the rack must be tangent to the pitch curves of both gears. 
See \reffig{TangentRack}.
The envelope of the rack during this motion forms the shape of the teeth that need to be cut out. 
See \reffig{Envelope}.
Since the pitch curve of the rack is not always perpendicular to the line connecting the two axles, the torque transmitted between the gears will not be constant during this motion. However this will produce teeth that have at least one point of contact at all times throughout the motion.

\subsection{Alien gears}
\label{Sec:AlienGears}

\begin{figure}[!h]
\centering
\subfloat[A strangely shaped gear.]{
\label{Fig:AlienStart}
\includegraphics[width = 0.4\textwidth]{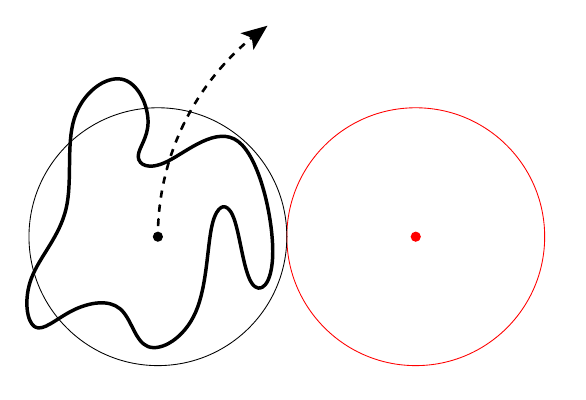}
}
\quad
\subfloat[Rolling the gear around.]{
\label{Fig:AlienRoll}
\includegraphics[width = 0.4\textwidth]{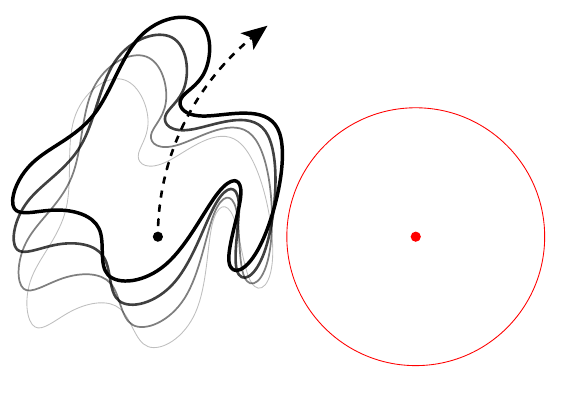}
}

\subfloat[The full sweep.]{
\label{Fig:AlienSweep}
\includegraphics[width = 0.6\textwidth]{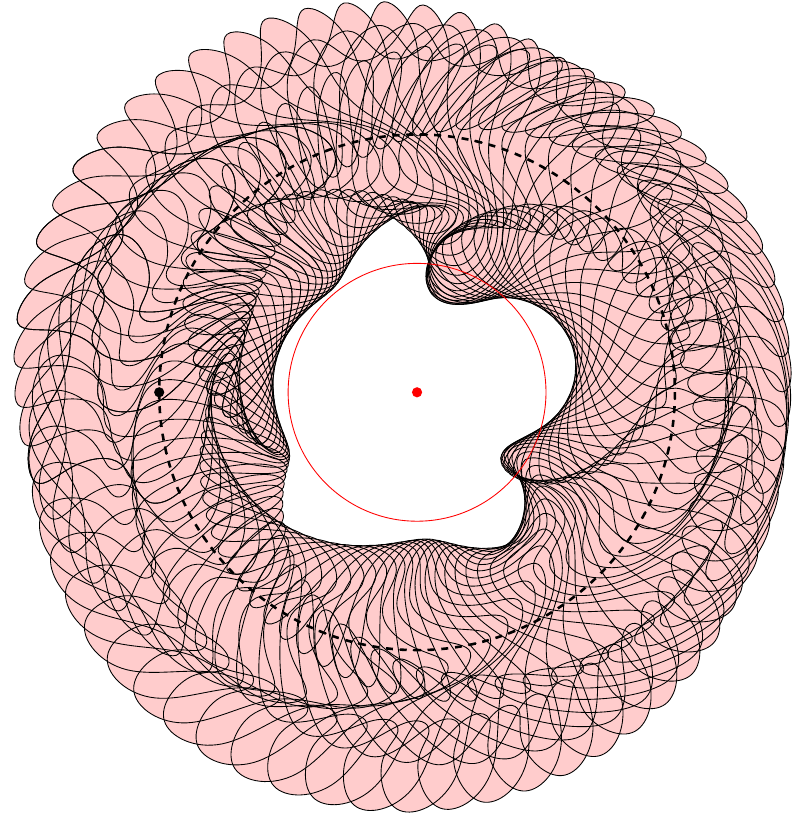}
}

\subfloat[]{
\label{Fig:RotateAlien1}
\includegraphics[width = 0.24\textwidth]{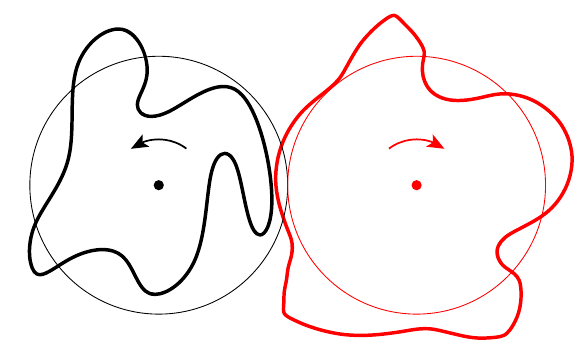}
}
\subfloat[]{
\label{Fig:RotateAlien2}
\includegraphics[width = 0.24\textwidth]{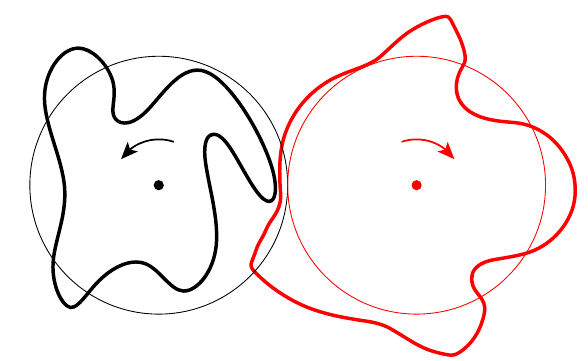}
}
\subfloat[]{
\label{Fig:RotateAlien3}
\includegraphics[width = 0.24\textwidth]{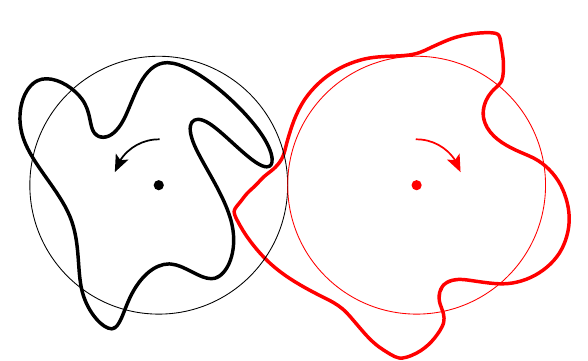}
}
\subfloat[]{
\label{Fig:RotateAlien4}
\includegraphics[width = 0.24\textwidth]{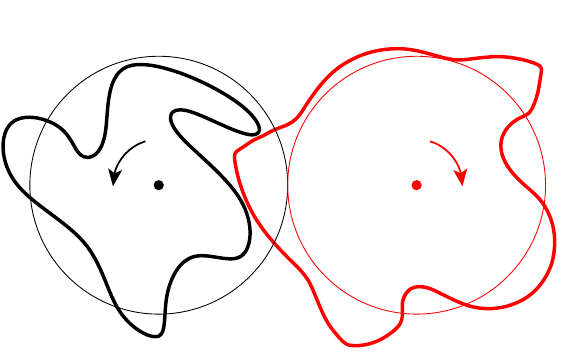}
}
\caption{Carving alien gears.}
\label{Fig:CarvingAlienGears}
\end{figure}

Alien gears are based on a pair of circular pitch curves with fixed radii $r_1$ and $r_2$. 
The (constant) rotation rates of the two gears are fixed by the ratio of their radii, but we take the shape of the ``teeth'' of the first gear to be some arbitrary shape, and must determine a corresponding shape for the second gear.
Conceptually, this can be achieved by carving the shape of the second gear from a large block of material by rotating the gears against each other at their constant rotation rates.
This is similar to the situation in \refsec{AcircularCutting}, but with the rack replaced by the first gear. 
In more detail, we begin with an arbitrary shape, $S_1$ and two gear centers, shown in \reffig{AlienStart}. 
We put ourselves in the frame of reference of the generated gear. 
The combined motion of the two rotations has $S_1$ roll around the second center, $c_2$, as shown in \reffig{AlienRoll}.
\reffig{AlienSweep} shows the full collection of positions for $S_1$ as it rolls around $c_2$. 
Assuming that the union of all of these positions has a hole in the middle containing $c_2$, we take $S_2$ to be the shape of the hole.

In Figures~\ref{Fig:RotateAlien1} to~\ref{Fig:RotateAlien4} we show the resulting gear shapes rotating against each other.
Depending on the choice of shape $S_1$, torque transfer will likely not be optimal -- indeed, the driving gear may not even force the driven gear to turn, and so the gears may slip against each other.

\begin{remark}
We made \reffig{CarvingAlienGears} using a Rhino-Grasshopper script by Arkuo Zheng. 
Rather than continuously sweeping $S_1$ around $c_2$, this uses a finite collection of positions for $S_1$.
It also smooths the shape of the hole in the center of \reffig{AlienSweep} to produce the gear shown in Figures~\ref{Fig:RotateAlien1} to~\ref{Fig:RotateAlien4}. 
For simple shapes and if we don't need much precision, this is good enough. 
Various approaches to computing continuous sweeps have been proposed.
See~\cite[Section~2]{MachchharSegermanElber} for an overview of work in this area.
\end{remark}

\section{Applications}

In this section, we describe some of our projects in 3D printed geared mechanisms.

\subsection{Odd numbers of gears}

Images of gears are often used in graphic design to give a sense of parts working together for a shared goal. 
In graphic design, the gears do not need to actually work -- in particular it is surprisingly common to see a design with an odd number of spur gears meshing in a loop. Adjacent spur gears rotate in opposite directions, so an odd number of these arranged in a loop is frozen in place.

\begin{figure}[!h]
\vspace{-20pt}
\centering
\subfloat[Triple gear by Saul Schleimer and the second author. Videos on this design:~\cite{TripleGearVideo, PoweredTripleGearVideo}]{
\label{Fig:TripleGear}
\includegraphics[width = 0.4\textwidth]{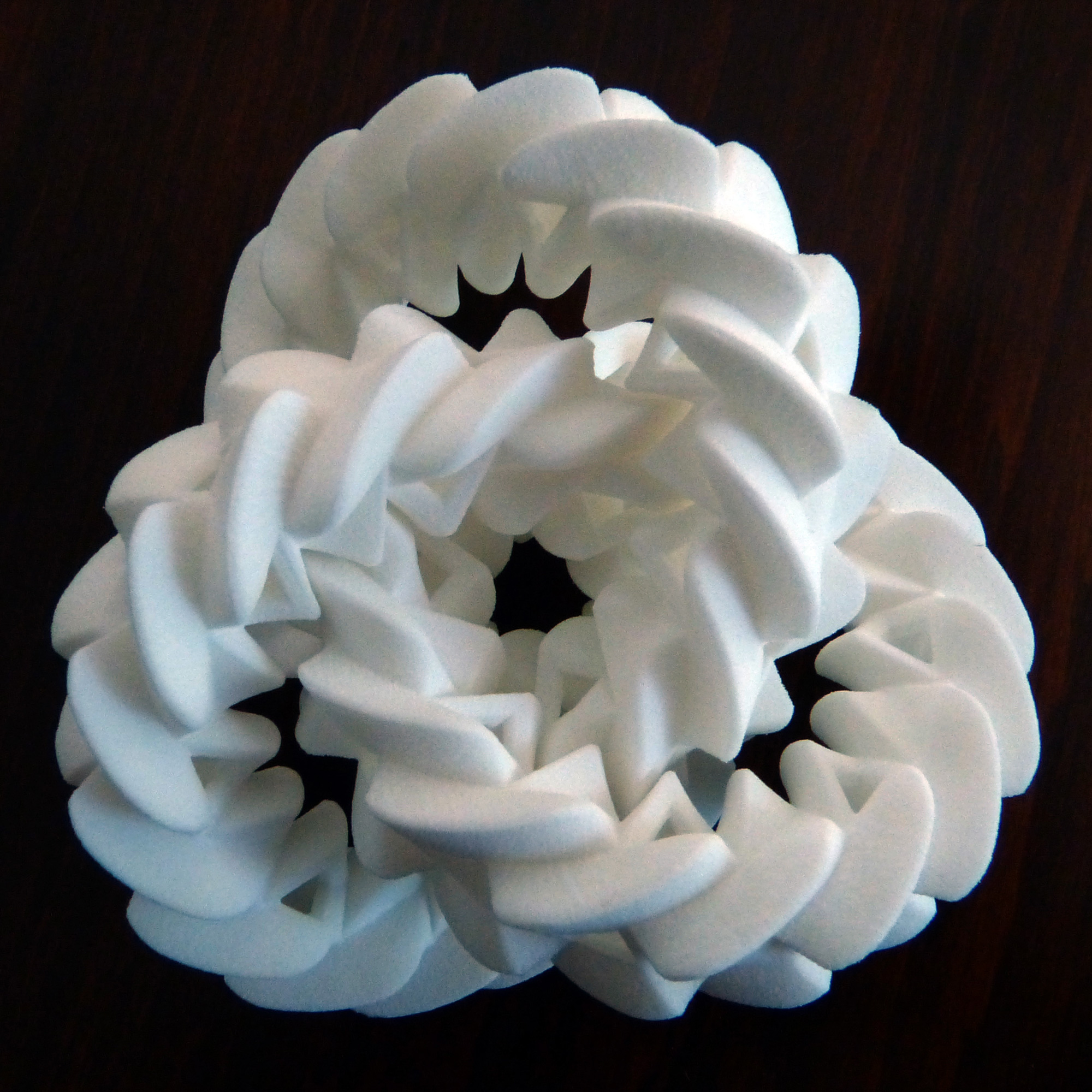}
}
\qquad
\subfloat[\emph{Triple helix} by Saul Schleimer and the second author. Video on this design:~\cite{TripleHelixVideo}]{
\includegraphics[height = 5cm]{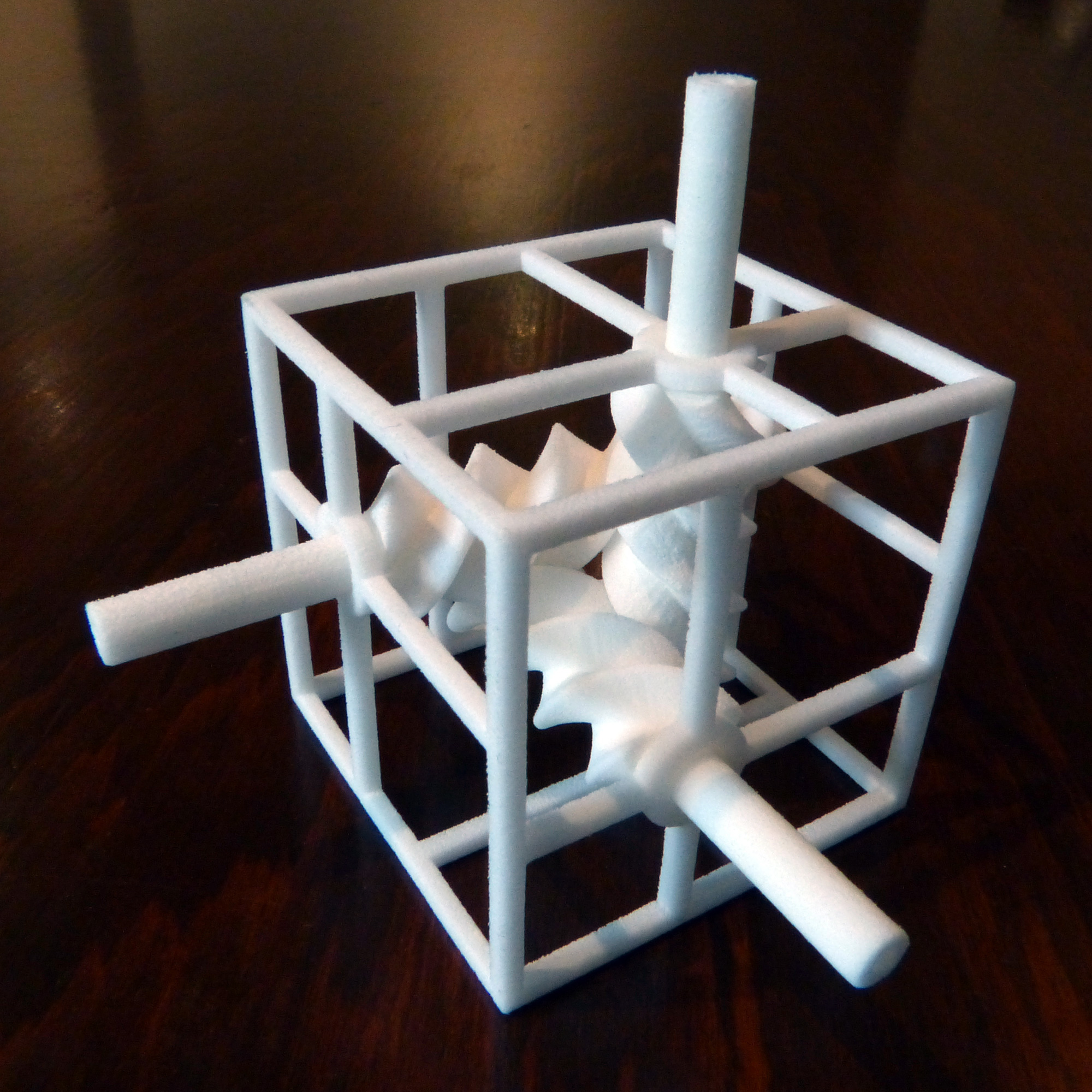}
\label{Fig:TripleHelix}
}
\caption{}
\label{Fig:Triples}
\end{figure}

An interesting challenge then is to construct a loop with an odd number of (necessarily) non-spur gears.
Our first solution was \emph{Triple gear}, see \reffig{TripleGear}. 
Here the three rings rotate around skew axes, so are skew gears. Unusually, some of the gearing surfaces face in towards the axle of the gear. See~\cite{TripleGear} for more details. A significantly simpler solution uses three helical gears, as shown in \reffig{TripleHelix}.

\subsection{Sliding racks}

\begin{figure}[!h]
\centering
\subfloat[\emph{Borromean racks}]{
\includegraphics[height = 4cm]{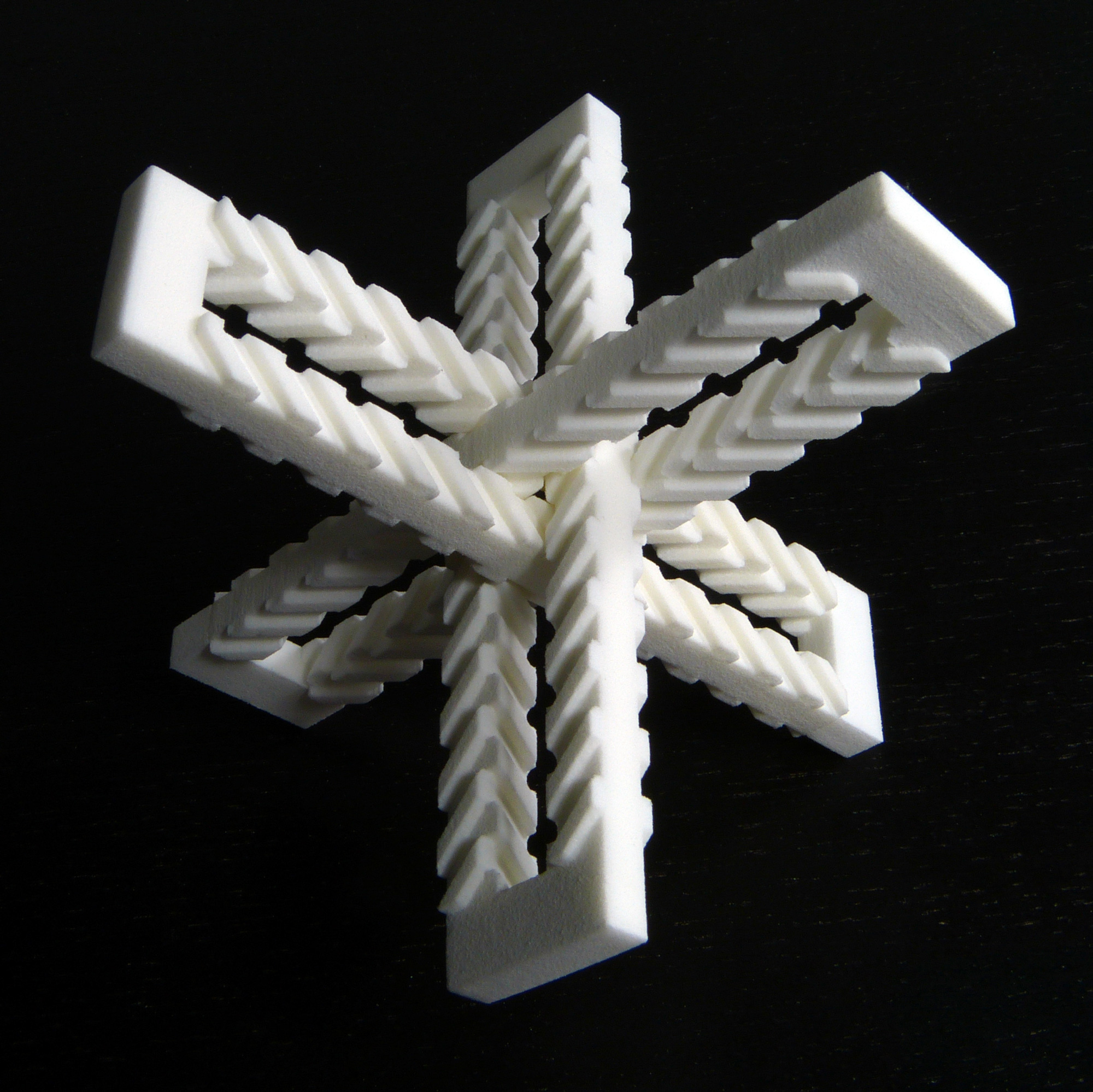}
}
\subfloat[\emph{Tetrahedral racks}]{
\includegraphics[height = 4cm]{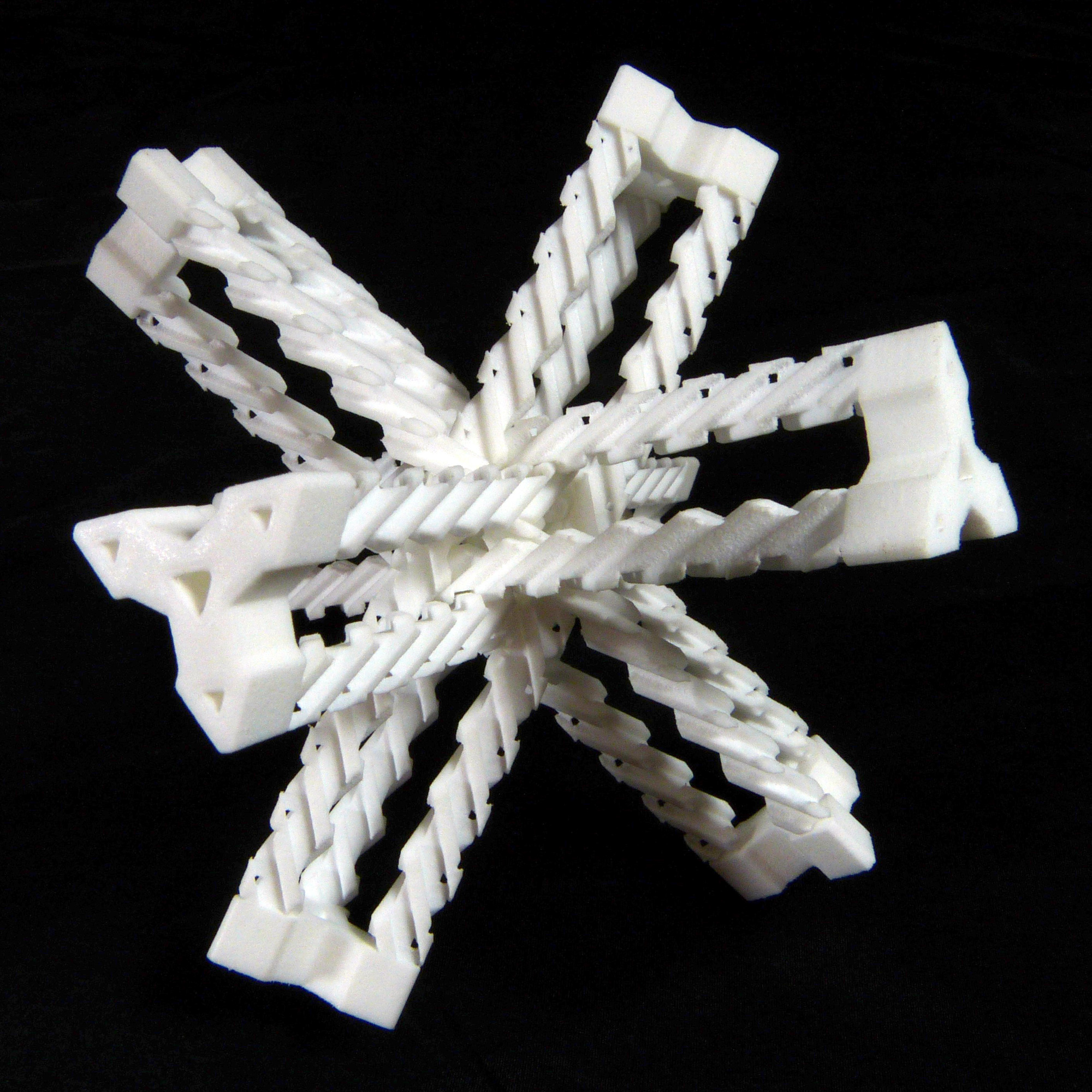}
}
\subfloat[\emph{Five axis racks}]{
\includegraphics[height = 4cm]{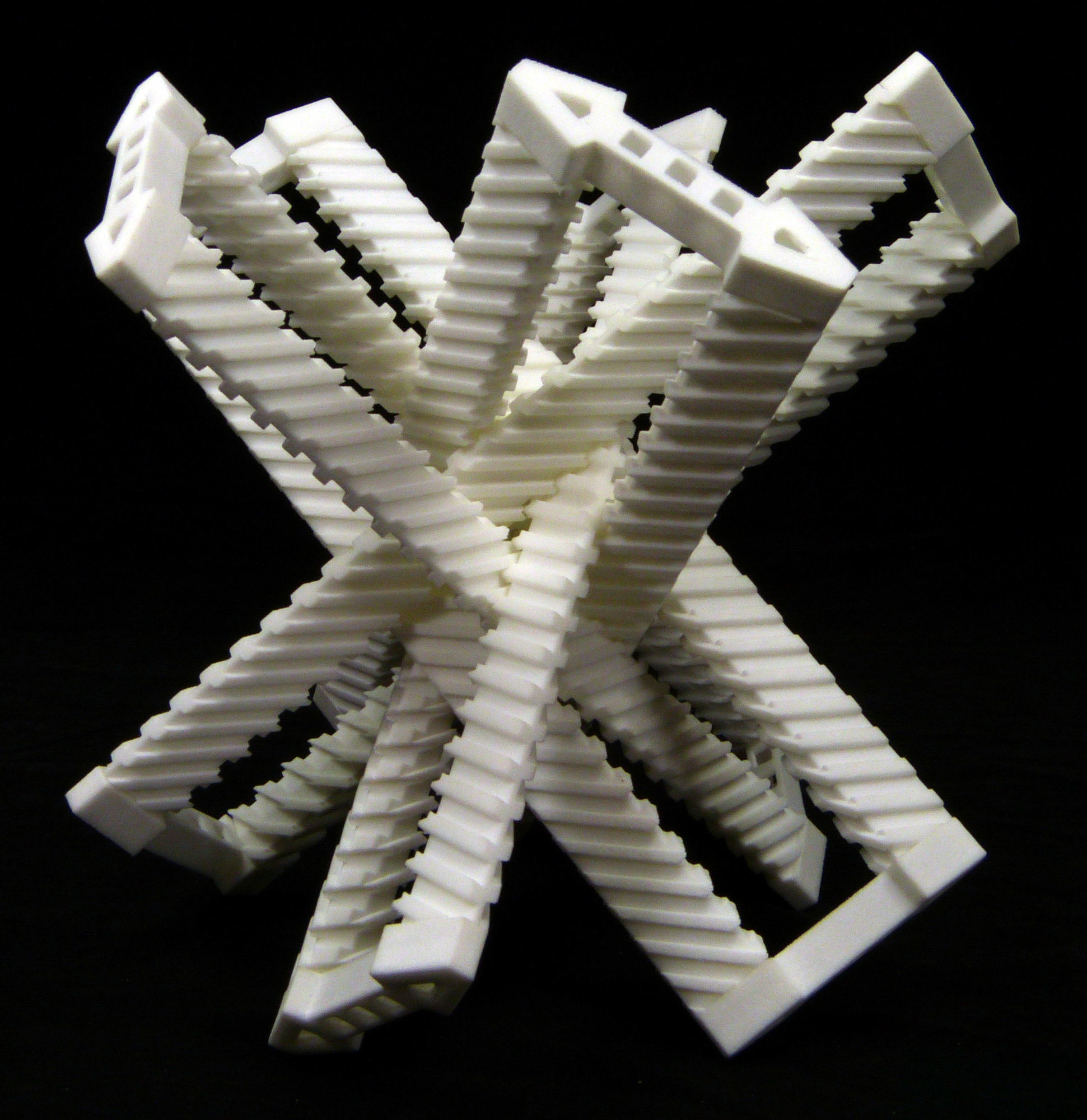}
}
\caption{Rack designs. The first two designs here are by Saul Schleimer and the second author, the third is by the authors. Videos on these designs:~\cite{BorromeanHairpinsVideo, TetrahedralRacksVideo, FiveAxisRacksVideo}}
\label{Fig:Racks}
\end{figure}

Racks allow a mechanism to convert the traditional rotational movement of a gear into translation. Mechanisms with gears and racks are usually designed to avoid parts sliding against each other. 
Sometimes however, sliding is the only way to achieve a desired outcome.
 \reffig{Racks} shows some of our mechanisms that use racks sliding against each other. The first of these, \emph{Borromean racks}, was also inspired by the odd-gear challenge.
  
\subsection{Gripping gears}

\begin{figure}[!h]
\centering
\includegraphics[width=0.8\textwidth]{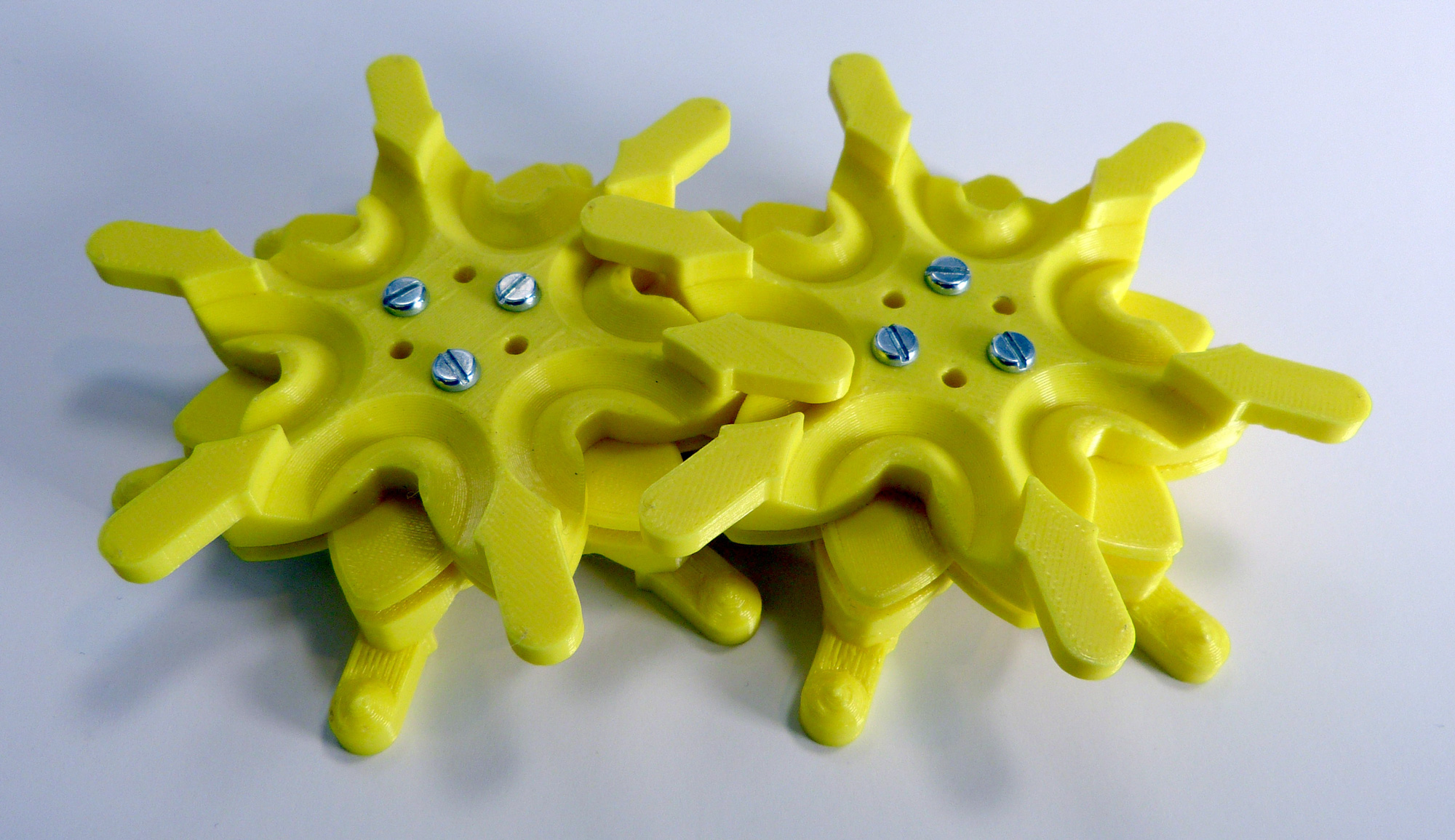}
\caption{\emph{Gripping gears}. Video on this design:~\cite{GrippingGearsVideo}}
\label{Fig:Gripping}
\end{figure}

Another challenge we set ourselves was to make a pair of gears that mesh and rotate together, be unable to separate, but have no axles attached to an outer frame. Our \emph{gripping gears} (joint work with Will Segerman) give a solution to this. See \reffig{Gripping}. 

Arms extending from each gear terminate in pegs which fit into grooves on the other gear. While a peg is in a groove, it restricts the movement of the two gears relative to each other. With enough such pegs and grooves, the two gears should be restricted to only be able to rotate around each other in the desired way. 

\begin{figure}[!h]
\centering
\includegraphics[width=0.8\textwidth]{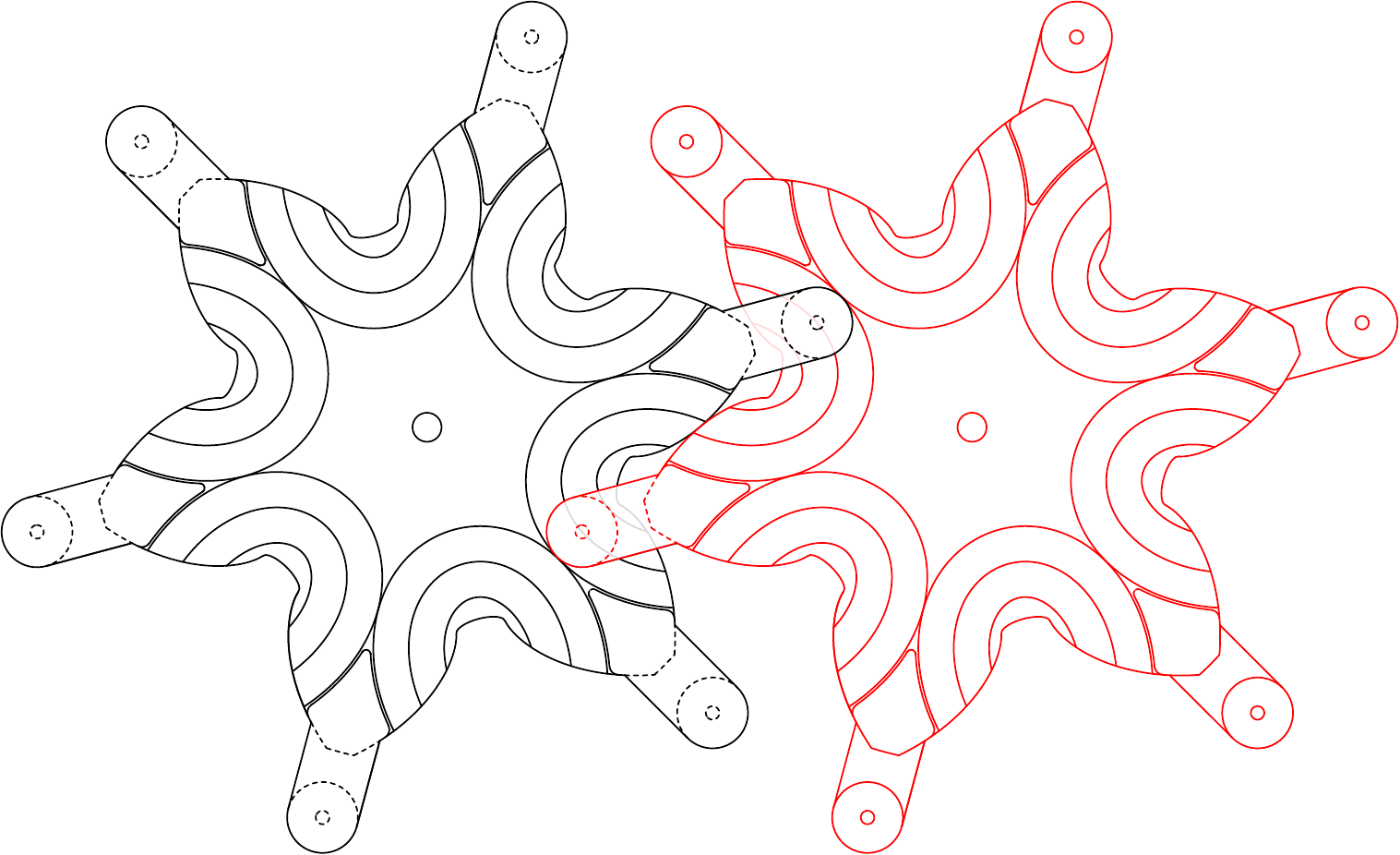}
\caption{Gripping gears schematic diagram.}
\label{Fig:GrippingSchematic}
\end{figure}

\begin{figure}[!h]
\centering
\subfloat[]{
\includegraphics[height=3cm]{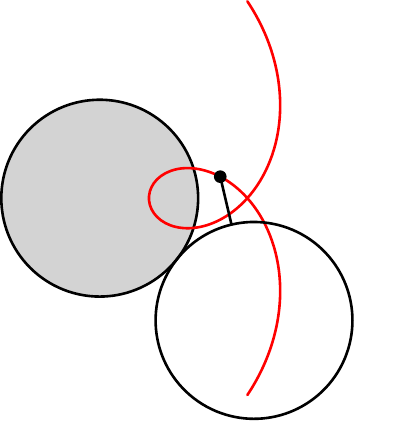}
}
\quad 
\subfloat[]{
\includegraphics[width=3cm]{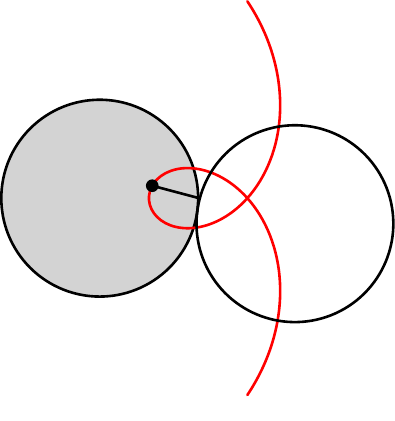}
}
\quad 
\subfloat[]{
\includegraphics[width=3cm]{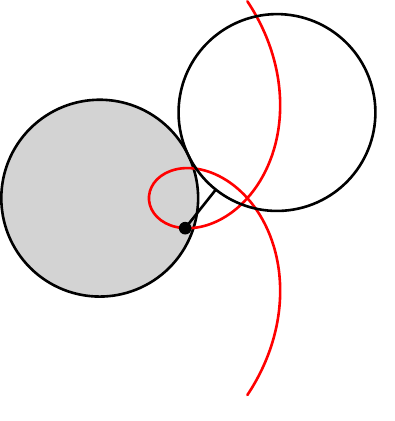}
}
\caption{Three frames of an animation: the white circle rolls around the outside of the grey circle. The end of an arm attached to the white circle traces out an epitrochoid in the frame of reference of the grey circle.}
\label{Fig:Epitrochoid}
\end{figure}

\reffig{GrippingSchematic} shows a schematic diagram of two pegs in their corresponding grooves. The shape of the groove is an \emph{epitrochoid}: the curve traced by a point attached to a circle as it rolls around a second, fixed circle. See \reffig{Epitrochoid}.

If trapped between two parallel planes, the two gears shown in \reffig{GrippingSchematic}, remain connected together for most of the desired motion. However, there is a configuration in which the two parts can be slid apart. In \reffig{GrippingSchematic}, if we slide the gear on the left further to the left and slightly up, then it can be disentangled from the gear on the right. Even worse, if the two gears are allowed to twist out of the plane then the pegs disengage and they can be untangled easily. To keep the parts together even when twisting is allowed, we make each part from two offset copies of the gear, sandwiched around two interleaving extra layers. The two copies protect each other's ``blind spots'' for planar disengagement, while the interleaving extra layers prevent twisting out of plane.
This then adds enough constraints to allow only the desired motion.

\begin{figure}[!h]
\centering
\subfloat[]{
\includegraphics[height = 5cm]{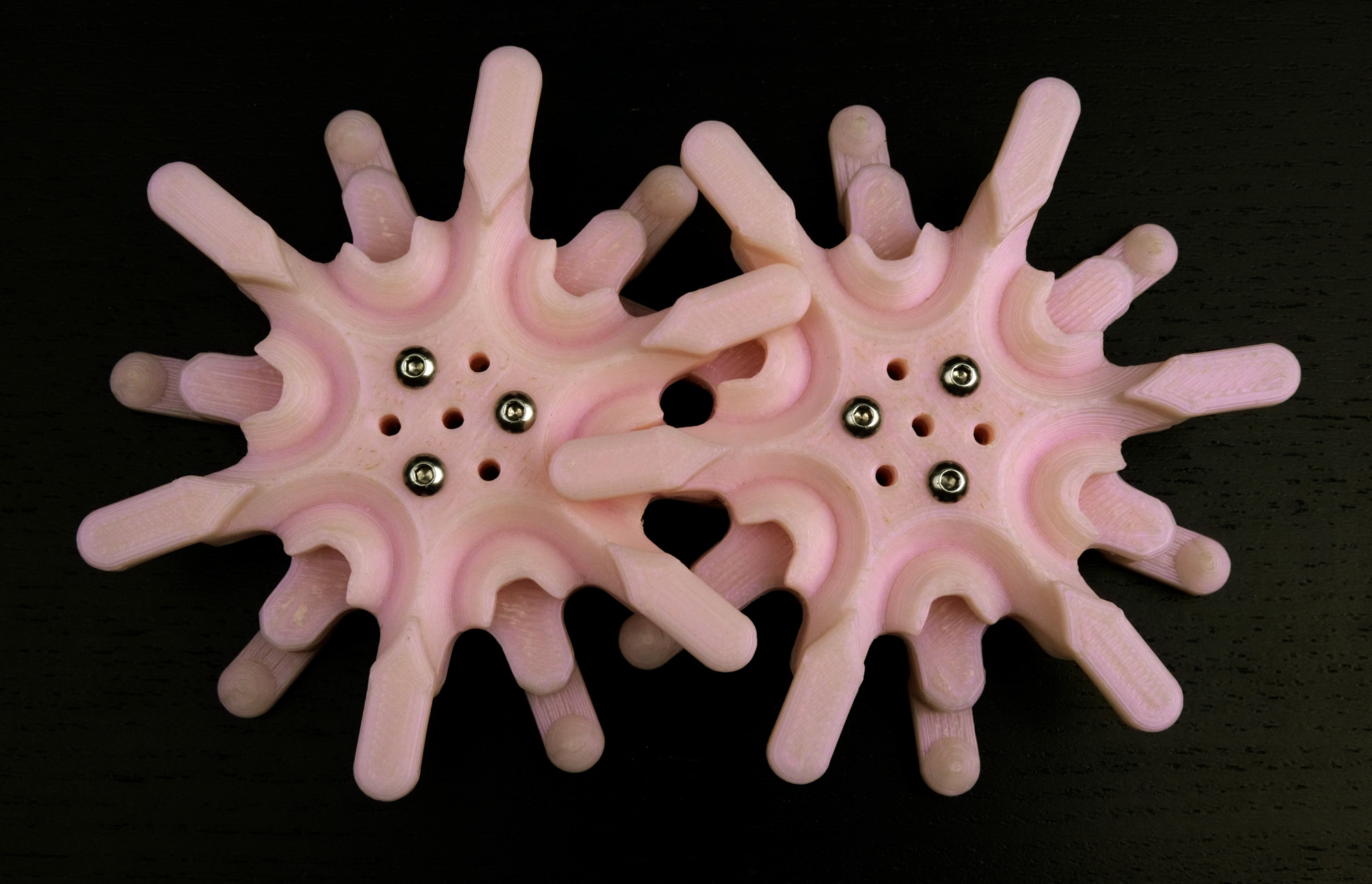}
}
\subfloat[]{
\includegraphics[height = 5cm]{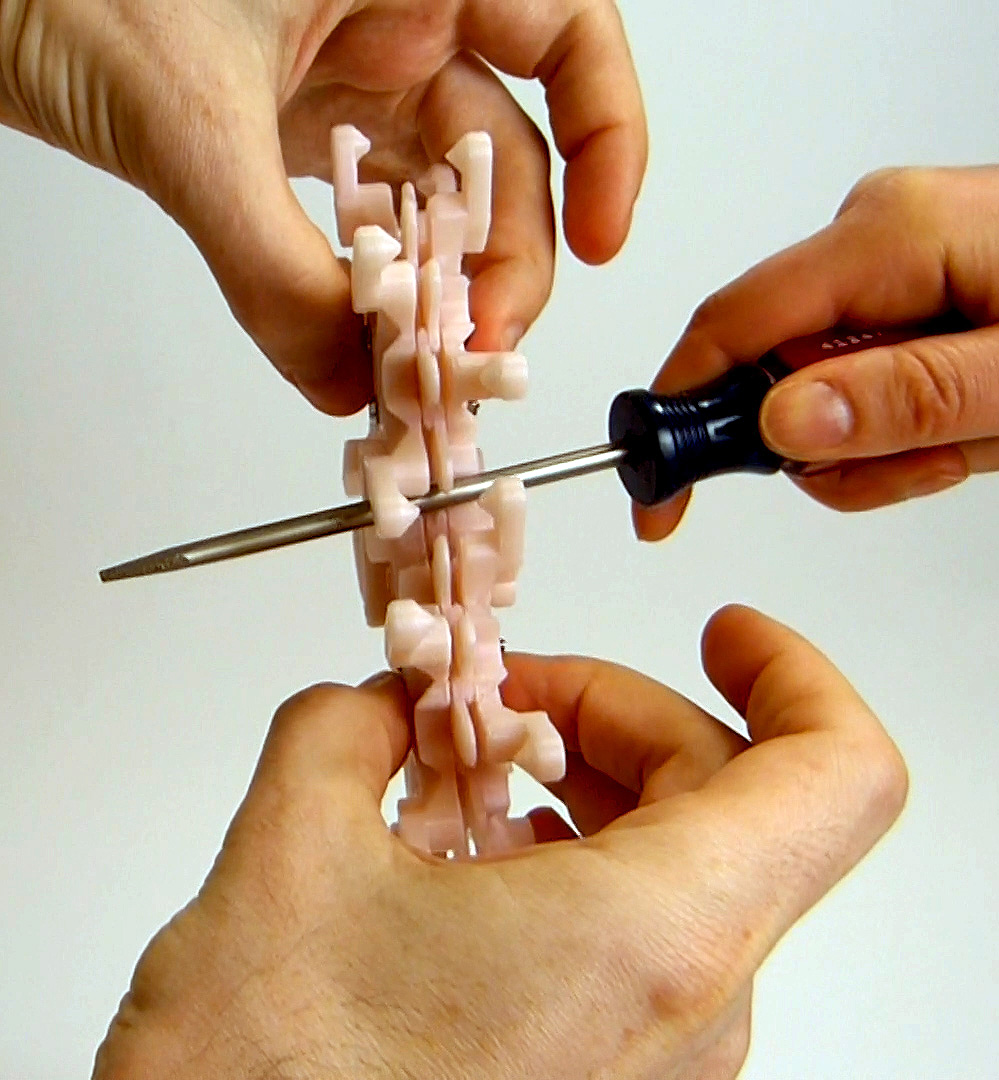}
\label{Fig:GrippingGearsWithHolesScrewdriver}
}
\caption{Gripping gears with holes. Video on this design:~\cite{GrippingGearsHolesVideo}}
\label{Fig:GrippingGearsWithHoles}
\end{figure}

In fact, there are enough constraints that we can remove some material and still have the mechanism work. In \reffig{GrippingGearsWithHoles} we have tunneled holes through some of the involute interfaces between the two gears. This makes the gear movement slightly less smooth but otherwise has little effect, and allows a solid rod (in \reffig{GrippingGearsWithHolesScrewdriver} a screwdriver) to pass between the two gears.

\subsection{Braiding gears}

One of our motivations in designing gripping gears was as a precursor to gears that would mesh and hold together for some part of their motion, but then disengage at some point. Our \emph{braiding gears} achieve this. 

\begin{figure}[!h]
\centering
\subfloat[]{
\includegraphics[width = \textwidth]{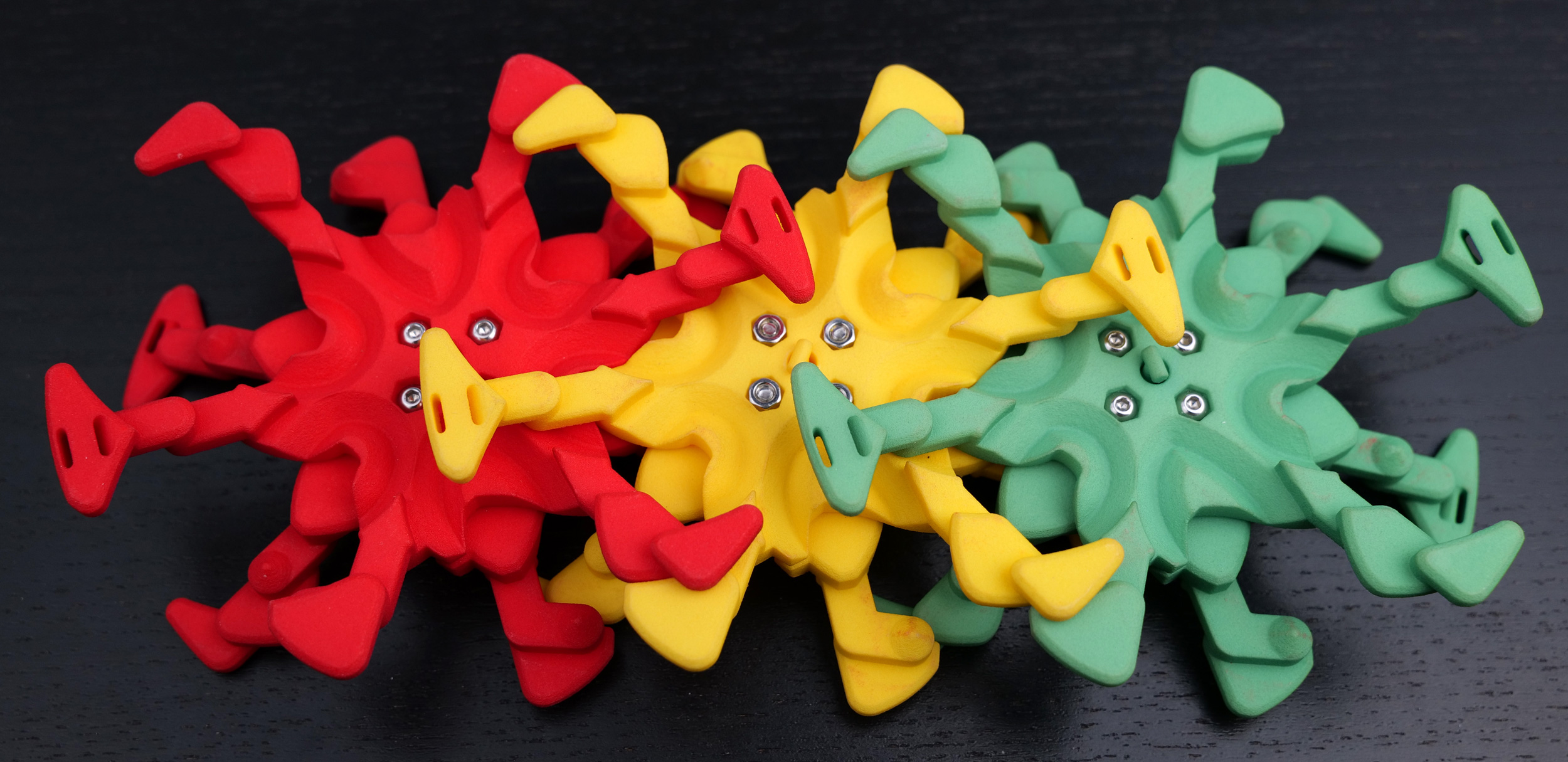}
\label{Fig:BraidingGears1}
}

\subfloat[]{
\includegraphics[height=3cm]{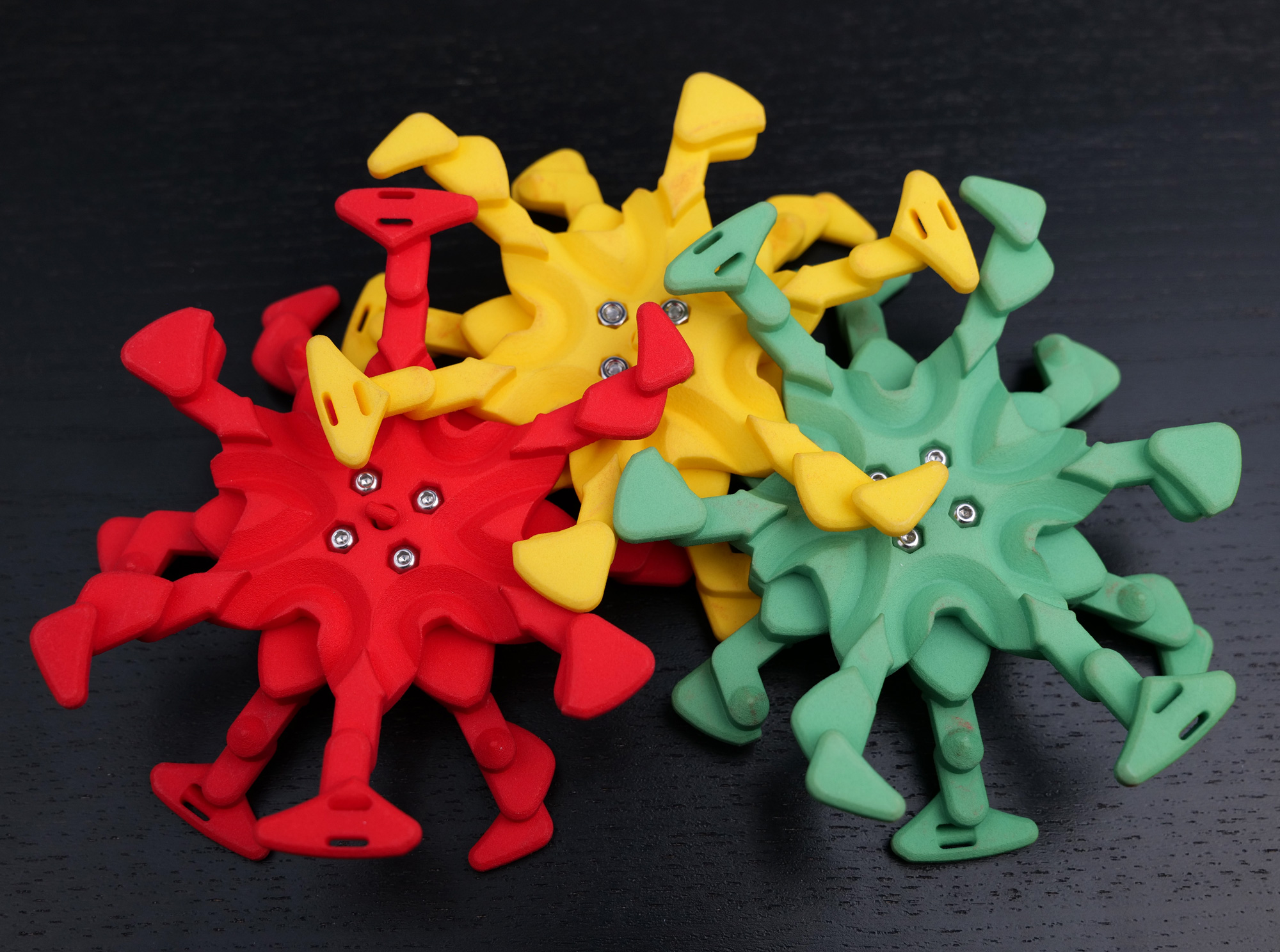}
\label{Fig:BraidingGears2}
}
\subfloat[]{
\includegraphics[height=3cm]{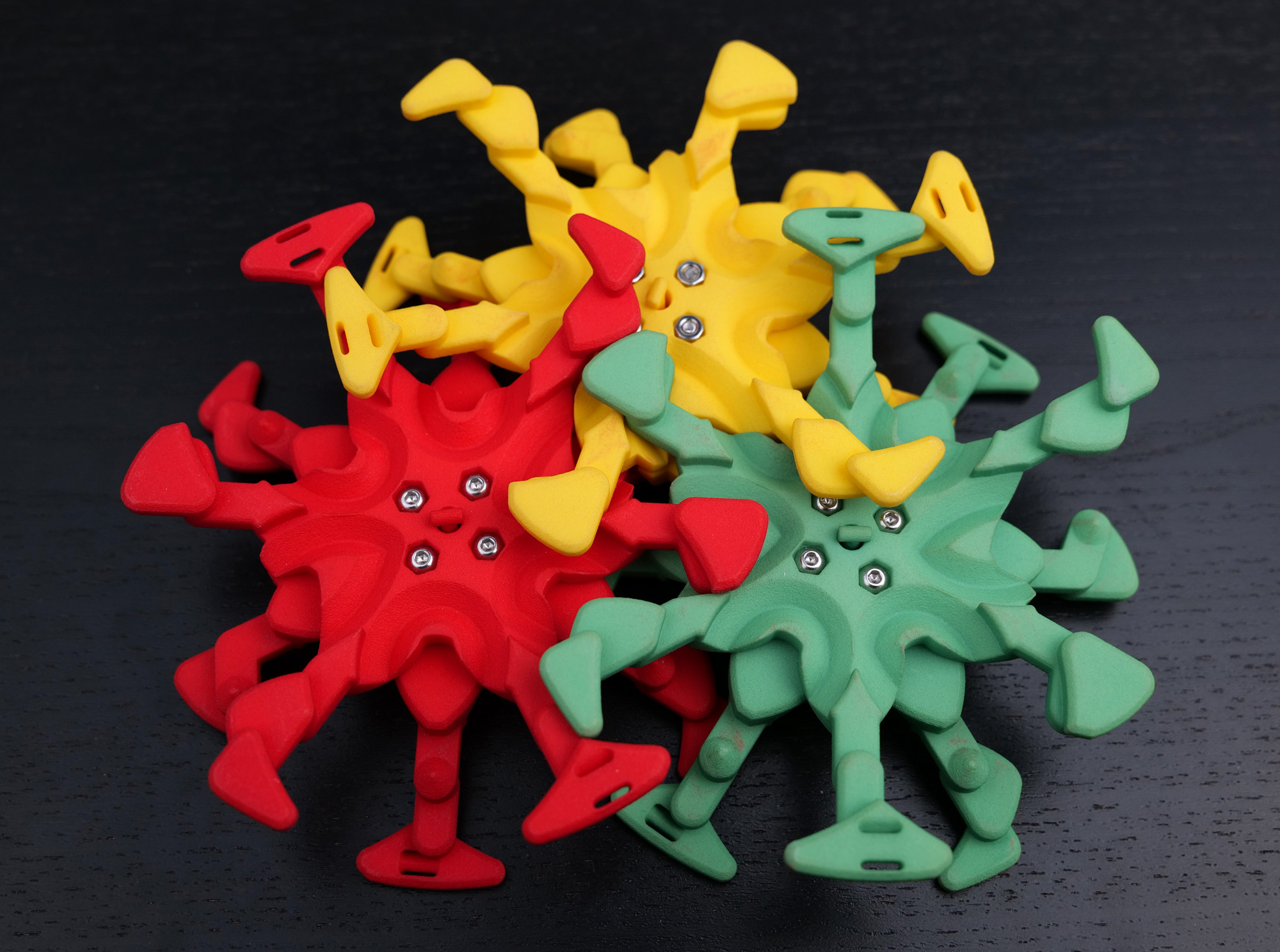}
\label{Fig:BraidingGears3}
}
\subfloat[]{
\includegraphics[height=3cm]{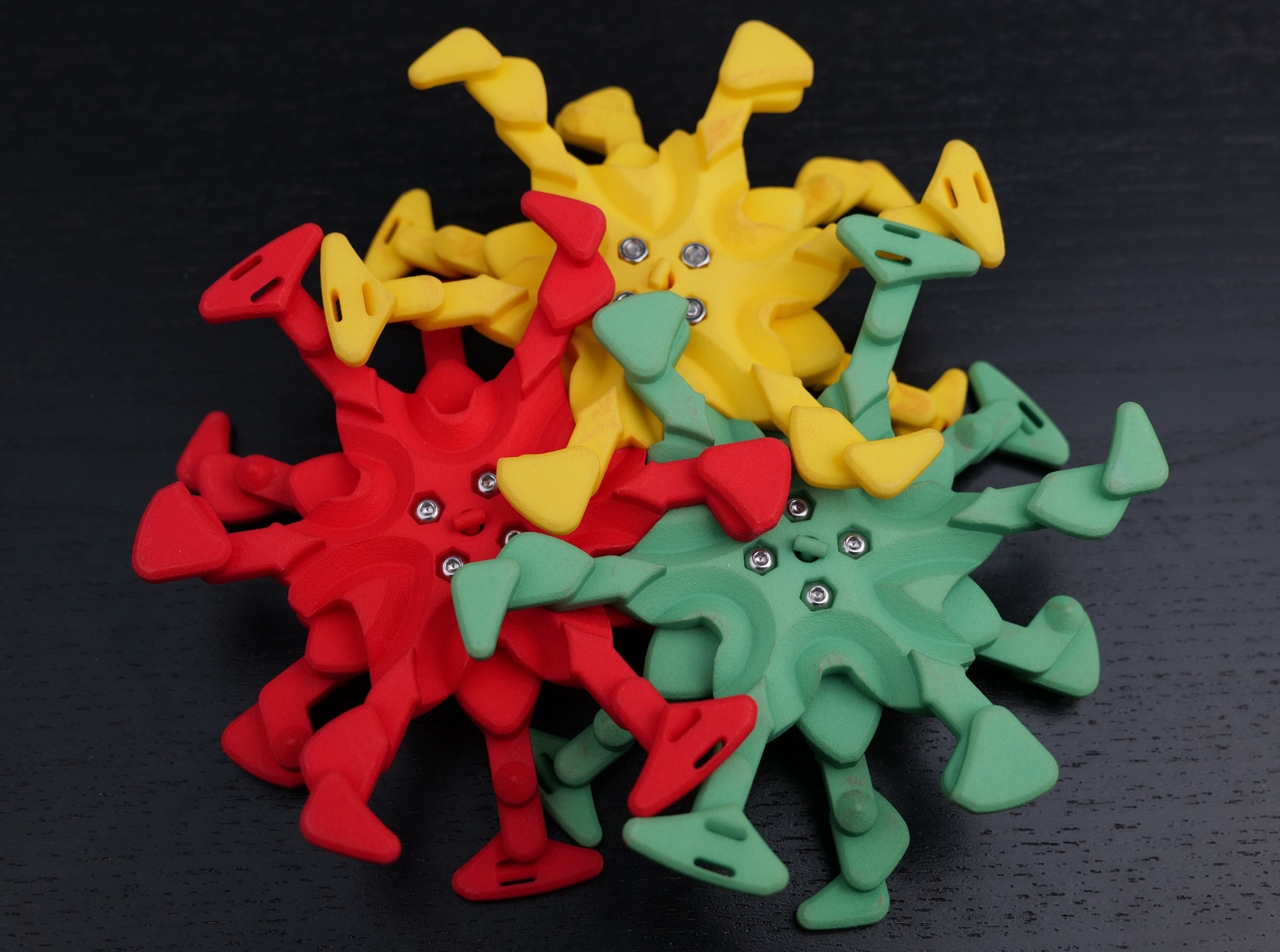}
\label{Fig:BraidingGears4}
}
\caption{\emph{Braiding gears}. Video on this design:~\cite{BraidingGearsVideo}}
\label{Fig:BraidingGears}
\end{figure}

In \reffig{BraidingGears1}, three identical gears are arranged in a line, with the outer two gears meshing only with the center gear. In \reffig{BraidingGears2} and \reffig{BraidingGears3} the outer two gears rotate around below the center gear, until in \reffig{BraidingGears4} they meet, and the configuration has three-fold rotational symmetry. (In fact, it has dihedral symmetry.) Thus, we may move away from this symmetric state in three different ways, each resulting in a different gear becoming the center gear in the line configuration. The video~\cite{BraidingGearsVideo} shows this motion.

\begin{figure}[!h]
\centering
\subfloat[]{
\includegraphics[width = 0.48\textwidth]{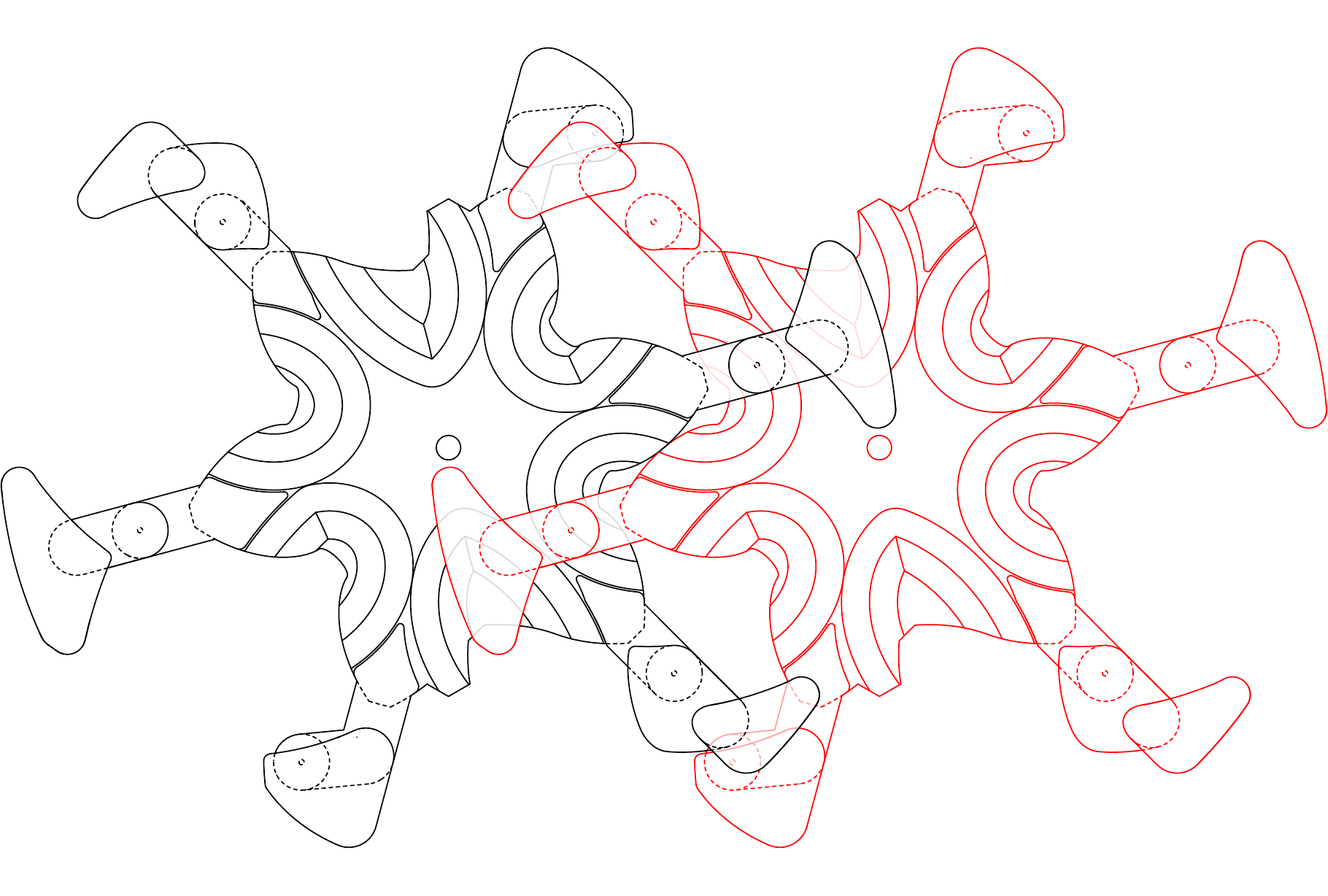}
\label{Fig:BraidingGearsSchematic1}
}
\subfloat[]{
\includegraphics[width = 0.48\textwidth]{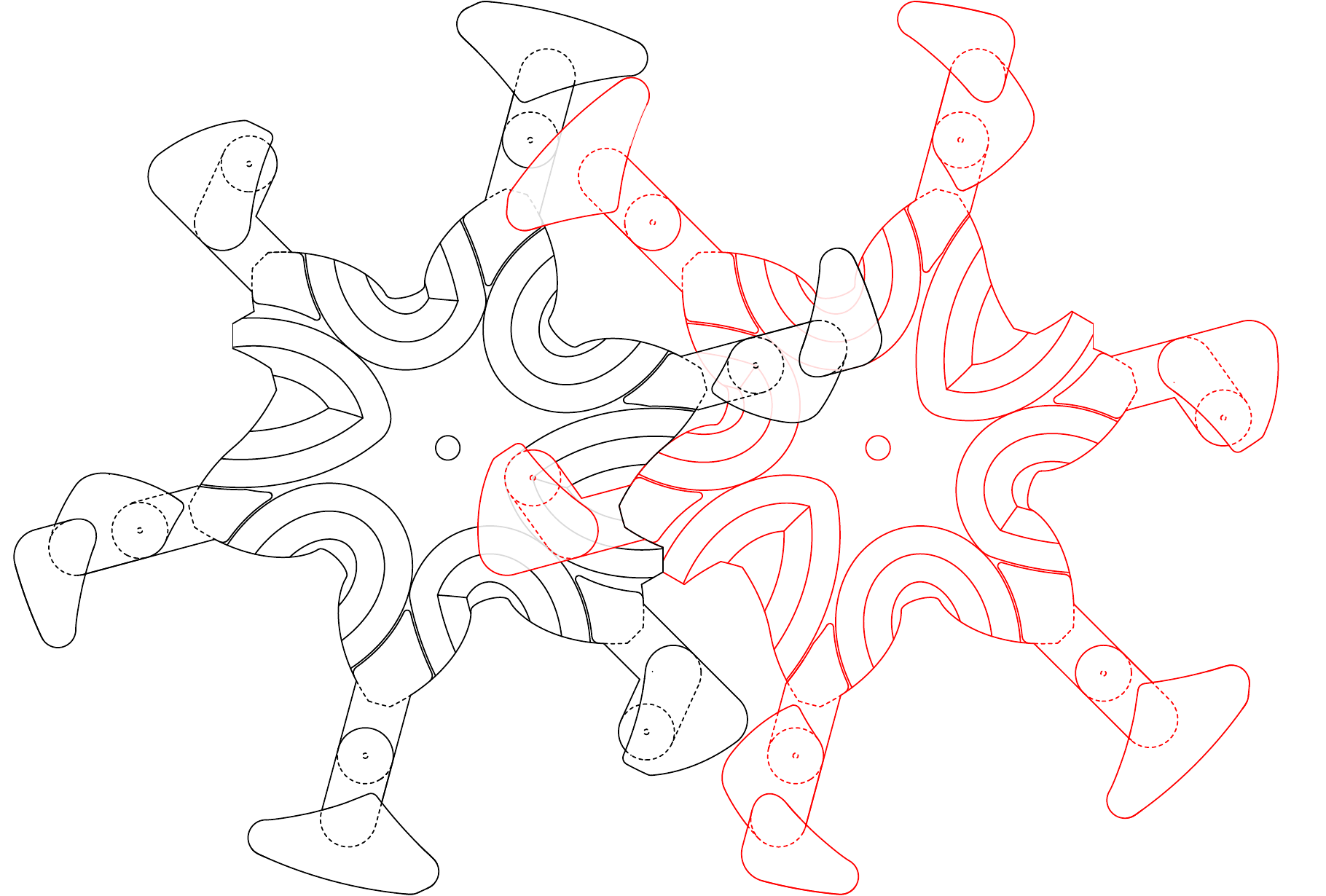}
\label{Fig:BraidingGearsSchematic2}
}
\caption{Braiding gears schematic diagram.}
\label{Fig:BraidingGearsSchematic}
\end{figure}

\reffig{BraidingGearsSchematic1} shows the upper side of two gears in the line configuration. Note that here they are translates of each other, and that their engaged pegs are in epicycloid grooves, just as for gripping gears. Rotating the left gear by $\pi/3$ anticlockwise and rotating the right gear by $\pi/3$ clockwise results in \reffig{BraidingGearsSchematic2}. The grooves in use here are made from two parts, meeting at a corner. The first part is a subarc of the standard epicycloid from the gripping gear design: a peg moves along this part as we rotate the gears from their positions in \reffig{BraidingGearsSchematic1} to those in \reffig{BraidingGearsSchematic2}. The second part is determined by the motion of the two gears that become outermost as we move from the three-fold symmetric configuration in \reffig{BraidingGears4} back to that in \reffig{BraidingGears1}: it is the path followed by a peg on one of those two gears as it moves relative to the other gear. 

\begin{figure}[!h]
\centering
\includegraphics[width=0.8\textwidth]{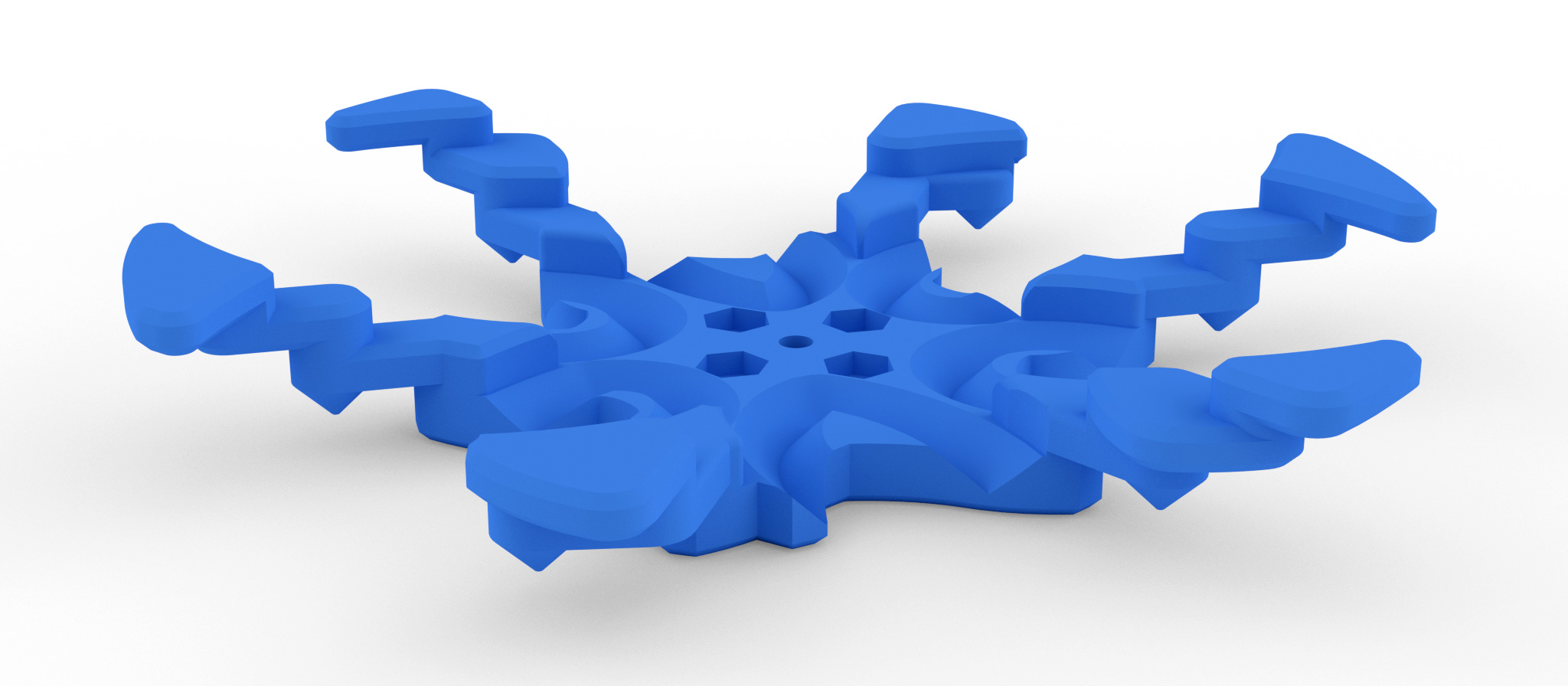}
\caption{Four levels.}
\label{Fig:BraidingGearsRender}
\end{figure}

These grooves and the corresponding pegs then allow the motion shown in \reffig{BraidingGears}. However, as shown in \reffig{BraidingGearsRender}, the design has parts extending beyond and above the grooves and pegs. Moving from bottom to top, there are a total of four ``levels'' to the design: level zero has the grooves in it, level one are the arms below which the pegs hang. Levels two and three are there to solve the following problem.

Since the configuration in \reffig{BraidingGears4} has three-fold symmetry, we can break the symmetry in three different ways. In particular, since the lower right gear could become the central gear, it is possible for the lower left gear to detach from the upper gear. If the lower right gear were not there to ``catch'' the lower left gear, then it would detach from the other two gears, and the mechanism would fall apart.

Thus, it is necessary that when in the line configuration, if one of the two outer gears rotates down then the other must also rotate down. In other words, if one outer gear rotates counterclockwise, then the other must rotate clockwise. This is the purpose of levels two and three -- on the two outer gears these levels interfere with each other, forcing the outer gears to stay in phase. 

We would have preferred if the two outer gears smoothly geared into each other. Unfortunately there is a lot of traffic moving around at the different levels as the mechanism moves, which restricts the size and shape of the upper level arms, and seems to preclude the possibility of perfectly smooth motion.

As in the gripping gears, the gears shown in \reffig{BraidingGearsSchematic} can be slid apart at particular points of their motion. Again, we solve this problem by bolting together two copies of the design, back to back. Here, the extra arms introduced on levels two and three have the side benefit of preventing twisting motions, and so the interleaving layers are not needed.

\bibliographystyle{hyperplain}
\bibliography{gears}
\end{document}

%% file: header_basic.tex
\usepackage{amsmath} 
\usepackage{amsthm} 
\usepackage{amssymb} 

\usepackage{microtype} 
\usepackage{pinlabel} 

\usepackage{MnSymbol} 

\usepackage[scaled=0.9]{sourcecodepro} 

\usepackage[hidelinks, pagebackref]{hyperref}

\makeatletter
\define@key{href}{font}{#1}
\makeatother
\usepackage{xpatch}
\newcommand\hrefdefaultfont{\ttfamily}
\xpatchcmd\href{\setkeys{href}{#1}}{\setkeys{href}{font=\hrefdefaultfont,#1}}{}{\fail}

\renewcommand*{\backref}[1]{}
\renewcommand*{\backrefalt}[4]{
  \ifcase #1 
  [No citations.]
  \or [#2]
  \else [#2]
  \fi }

\let\originalleft\left
\let\originalright\right
\renewcommand{\left}{\mathopen{}\mathclose\bgroup\originalleft}
\renewcommand{\right}{\aftergroup\egroup\originalright}

%% file: header_subtle.tex
\theoremstyle{plain}
\newtheorem{XXXtheoremQED}[equation]{Theorem} 
  {\pushQED{\qed}\begin{XXXtheoremQED}}
  {\popQED\end{XXXtheoremQED}}

\newcommand{\fakeenv}{} 

{ 
 \renewcommand{\fakeenv}{#2} 
 \theoremstyle{plain} 
 \newtheorem*{\fakeenv}{#1~\ref{#2}} 
 \begin{\fakeenv}
}
{
 \end{\fakeenv}
}

\newenvironment{restated}[2]  
{ 
 \renewcommand{\fakeenv}{#2} 
 \theoremstyle{definition} 
 \newtheorem*{\fakeenv}{#1~\ref{#2}} 
 \begin{\fakeenv}
}
{
 \end{\fakeenv}
}








%% file: gear_design.bbl
\begin{thebibliography}{10}

\bibitem{Bencsik}
Gergely Bencsik.
\newblock Involute gears explained, July 2022.
\newblock \url{https://www.youtube.com/watch?v=nrsCoQN6V4M}.

\bibitem{acirc_involute}
B{\'a}lint Laczik.
\newblock Involute profile of non-circular gears.
\newblock Available from
  \url{http://manuals.chudov.com/Non-Circular-Gears.pdf}.

\bibitem{MachchharSegermanElber}
Jinesh Machchhar, Henry Segerman, and Gershon Elber.
\newblock Conjugate shape simplification via precise algebraic planar sweeps
  toward gear design.
\newblock {\em Computers \& Graphics}, 90:1--10, 2020.

\bibitem{FiveAxisRacksVideo}
Elisabetta~A. Matsumoto and Henry Segerman.
\newblock Five axis racks, September 2018.
\newblock \url{https://www.youtube.com/watch?v=dkguSyeQXjc}.

\bibitem{GearedJitterbugVideo}
Elisabetta~A. Matsumoto and Henry Segerman.
\newblock Geared jitterbug, October 2018.
\newblock \url{https://www.youtube.com/watch?v=kwERR5flAOU}.

\bibitem{GearedCuboctahedralJitterbugVideo}
Elisabetta~A. Matsumoto and Henry Segerman.
\newblock Geared cuboctahedral jitterbug, March 2019.
\newblock \url{https://www.youtube.com/watch?v=eZFgHnvtPUo}.

\bibitem{GearedJitterbugs}
Elisabetta~A. Matsumoto and Henry Segerman.
\newblock Geared jitterbugs.
\newblock In Susan Goldstine, Douglas McKenna, and Krist\'{o}f Fenyvesi,
  editors, {\em Proceedings of Bridges 2019: Mathematics, Art, Music,
  Architecture, Education, Culture}, pages 399--402, Phoenix, Arizona, 2019.
  Tessellations Publishing.

\bibitem{GrippingGearsHolesVideo}
Elisabetta~A. Matsumoto and Henry Segerman.
\newblock Gripping gears with pass-through holes, October 2019.
\newblock \url{https://www.youtube.com/watch?v=RBZG8M8_a8Y}.

\bibitem{BevelGearsVideo}
Elisabetta~A. Matsumoto and Henry Segerman.
\newblock Gear cube and {B}rain gear, December 2022.
\newblock \url{https://www.youtube.com/watch?v=l4eL1MvEbN0}.

\bibitem{GrippingGearsVideo}
Elisabetta~A. Matsumoto, Henry Segerman, and Will Segerman.
\newblock Gripping gears, September 2019.
\newblock \url{https://www.youtube.com/watch?v=ENFXnNtd3xU}.

\bibitem{ParkLee}
Noh~Gill Park and Hyoung~Woo Lee.
\newblock The spherical involute bevel gear: its geometry, kinematic behavior
  and standardization.
\newblock {\em Journal of Mechanical Science and Technology}, 25(4):1023--1034,
  2011.

\bibitem{TripleGearVideo}
Saul Schleimer and Henry Segerman.
\newblock Triple gear, December 2012.
\newblock \url{https://www.youtube.com/watch?v=I9IBQVHFeQs}.

\bibitem{BorromeanHairpinsVideo}
Saul Schleimer and Henry Segerman.
\newblock Borromean hairpins, September 2013.
\newblock \url{https://www.youtube.com/watch?v=WQ9ptuUxfk4}.

\bibitem{PoweredTripleGearVideo}
Saul Schleimer and Henry Segerman.
\newblock Powered triple gear, April 2013.
\newblock \url{https://www.youtube.com/watch?v=QhXjevOY_uk}.

\bibitem{TripleGear}
Saul Schleimer and Henry Segerman.
\newblock Triple gear.
\newblock In George~W. Hart and Reza Sarhangi, editors, {\em Proceedings of
  Bridges 2013: Mathematics, Music, Art, Architecture, Culture}, pages
  353--360, Phoenix, Arizona, 2013. Tessellations Publishing.

\bibitem{TripleHelixVideo}
Saul Schleimer and Henry Segerman.
\newblock Triple helix, March 2013.
\newblock \url{https://www.youtube.com/watch?v=74ygvXzmrgk}.

\bibitem{TetrahedralRacksVideo}
Saul Schleimer and Henry Segerman.
\newblock Tetrahedral racks, August 2016.
\newblock \url{https://www.youtube.com/watch?v=sDJk-9QxFro}.

\bibitem{BraidingGearsVideo}
Henry Segerman.
\newblock Braiding gears, November 2019.
\newblock \url{https://www.youtube.com/watch?v=Lh7yAbw0H24}.

\end{thebibliography}
